\documentclass[a4paper,12pt]{article}
\usepackage[all]{xy}
\usepackage[T1]{fontenc}
\usepackage[dvips]{graphicx}
\usepackage{amsmath,amssymb,amsfonts,amsthm,latexsym,mathrsfs,textcomp,verbatim}
\xyoption{web}

\title{De Morgan classifying toposes}
\author{{Olivia Caramello} \vspace{3 mm}\\ {\small DPMMS, University of Cambridge,}\\{\small Wilberforce Road, Cambridge CB3 0WB, UK}\\{\small O.Caramello@dpmms.cam.ac.uk}}
\date{August 11, 2008}
\begin{document}
\bgroup           
\let\footnoterule\relax  
\maketitle
\flushleft  
\begin{abstract}
We present a general method for deciding whether a Grothendieck topos satisfies De Morgan's law (resp. the law of excluded middle) or not; applications to the theory of classifying toposes follow. Specifically, we obtain a syntactic characterization of the class of geometric theories whose classifying toposes satisfy De Morgan's law (resp. are Boolean), as well as model-theoretic criteria for theories whose classifying toposes arise as localizations of a given presheaf topos.
\end{abstract} 
\egroup 
\flushleft
\vspace{5 mm}


\def\Monthnameof#1{\ifcase#1\or
   January\or February\or March\or April\or May\or June\or
   July\or August\or September\or October\or November\or December\fi}
\def\today{\number\day~\Monthnameof\month~\number\year}

%
%
%
\def\pushright#1{{
   \parfillskip=0pt            
   \widowpenalty=10000         
   \displaywidowpenalty=10000  
   \finalhyphendemerits=0      
  %
   \leavevmode                 
   \unskip                     
   \nobreak                    
   \hfil                       
   \penalty50                  
   \hskip.2em                  
   \null                       
   \hfill                      
   {#1}                        
  %
   \par}}                      

\def\qed{\pushright{$\square$}\penalty-700 \smallskip}

\newtheorem{theorem}{Theorem}[section]

\newtheorem{proposition}[theorem]{Proposition}

\newtheorem{scholium}[theorem]{Scholium}

\newtheorem{lemma}[theorem]{Lemma}

\newtheorem{corollary}[theorem]{Corollary}

\newtheorem{conjecture}[theorem]{Conjecture}

\newenvironment{proofs}%
 {\begin{trivlist}\item[]{\bf Proof }}%
 {\qed\end{trivlist}}

  \newtheorem{rmk}[theorem]{Remark}
\newenvironment{remark}{\begin{rmk}\em}{\end{rmk}}

  \newtheorem{rmks}[theorem]{Remarks}
\newenvironment{remarks}{\begin{rmks}\em}{\end{rmks}}

  \newtheorem{defn}[theorem]{Definition}
\newenvironment{definition}{\begin{defn}\em}{\end{defn}}

  \newtheorem{eg}[theorem]{Example}
\newenvironment{example}{\begin{eg}\em}{\end{eg}}

  \newtheorem{egs}[theorem]{Examples}
\newenvironment{examples}{\begin{egs}\em}{\end{egs}}


\mathcode`\<="4268  
\mathcode`\>="5269  
\mathcode`\.="313A  
\mathchardef\semicolon="603B 
\mathchardef\gt="313E
\mathchardef\lt="313C

\newcommand{\app}
 {{\sf app}}

\newcommand{\Ass}
 {{\bf Ass}}

\newcommand{\ASS}
 {{\mathbb A}{\sf ss}}

\newcommand{\Bb}
{\mathbb}

\newcommand{\biimp}
 {\!\Leftrightarrow\!}

\newcommand{\bim}
 {\rightarrowtail\kern-1em\twoheadrightarrow}

\newcommand{\bjg}
 {\mathrel{{\dashv}\,{\vdash}}}

\newcommand{\bstp}[3]
 {\mbox{$#1\! : #2 \bim #3$}}

\newcommand{\cat}
 {\!\mbox{\t{\ }}}

\newcommand{\cinf}
 {C^{\infty}}

\newcommand{\cinfrg}
 {\cinf\hy{\bf Rng}}

\newcommand{\cocomma}[2]
 {\mbox{$(#1\!\uparrow\!#2)$}}

\newcommand{\cod}
 {{\rm cod}}

\newcommand{\comma}[2]
 {\mbox{$(#1\!\downarrow\!#2)$}}

\newcommand{\comp}
 {\circ}

\newcommand{\cons}
 {{\sf cons}}

\newcommand{\Cont}
 {{\bf Cont}}

\newcommand{\ContE}
 {{\bf Cont}_{\cal E}}

\newcommand{\ContS}
 {{\bf Cont}_{\cal S}}

\newcommand{\cover}
 {-\!\!\triangleright\,}

\newcommand{\cstp}[3]
 {\mbox{$#1\! : #2 \cover #3$}}

\newcommand{\Dec}
 {{\rm Dec}}

\newcommand{\DEC}
 {{\mathbb D}{\sf ec}}

\newcommand{\den}[1]
 {[\![#1]\!]}

\newcommand{\Desc}
 {{\bf Desc}}

\newcommand{\dom}
 {{\rm dom}}

\newcommand{\Eff}
 {{\bf Eff}}

\newcommand{\EFF}
 {{\mathbb E}{\sf ff}}

\newcommand{\empstg}
 {[\,]}

\newcommand{\epi}
 {\twoheadrightarrow}

\newcommand{\estp}[3]
 {\mbox{$#1 \! : #2 \epi #3$}}

\newcommand{\ev}
 {{\rm ev}}

\newcommand{\Ext}
 {{\rm Ext}}

\newcommand{\fr}
 {\sf}

\newcommand{\fst}
 {{\sf fst}}

\newcommand{\fun}[2]
 {\mbox{$[#1\!\to\!#2]$}}

\newcommand{\funs}[2]
 {[#1\!\to\!#2]}

\newcommand{\Gl}
 {{\bf Gl}}

\newcommand{\hash}
 {\,\#\,}

\newcommand{\hy}
 {\mbox{-}}

\newcommand{\im}
 {{\rm im}}

\newcommand{\imp}
 {\!\Rightarrow\!}

\newcommand{\Ind}[1]
 {{\rm Ind}\hy #1}

\newcommand{\iten}[1]
{\item[{\rm (#1)}]}

\newcommand{\iter}
 {{\sf iter}}

\newcommand{\Kalg}
 {K\hy{\bf Alg}}

\newcommand{\llim}
 {{\mbox{$\lower.95ex\hbox{{\rm lim}}$}\atop{\scriptstyle
{\leftarrow}}}{}}

\newcommand{\llimd}
 {\lower0.37ex\hbox{$\pile{\lim \\ {\scriptstyle
\leftarrow}}$}{}}

\newcommand{\Mf}
 {{\bf Mf}}

\newcommand{\Mod}
 {{\bf Mod}}

\newcommand{\MOD}
{{\mathbb M}{\sf od}}

\newcommand{\mono}
 {\rightarrowtail}

\newcommand{\mor}
 {{\rm mor}}

\newcommand{\mstp}[3]
 {\mbox{$#1\! : #2 \mono #3$}}

\newcommand{\Mu}
 {{\rm M}}

\newcommand{\name}[1]
 {\mbox{$\ulcorner #1 \urcorner$}}

\newcommand{\names}[1]
 {\mbox{$\ulcorner$} #1 \mbox{$\urcorner$}}

\newcommand{\nml}
 {\triangleleft}

\newcommand{\ob}
 {{\rm ob}}

\newcommand{\op}
 {^{\rm op}}

\newcommand{\pepi}
 {\rightharpoondown\kern-0.9em\rightharpoondown}

\newcommand{\pmap}
 {\rightharpoondown}

\newcommand{\Pos}
 {{\bf Pos}}

\newcommand{\prarr}
 {\rightrightarrows}

\newcommand{\princfil}[1]
 {\mbox{$\uparrow\!(#1)$}}

\newcommand{\princid}[1]
 {\mbox{$\downarrow\!(#1)$}}

\newcommand{\prstp}[3]
 {\mbox{$#1\! : #2 \prarr #3$}}

\newcommand{\pstp}[3]
 {\mbox{$#1\! : #2 \pmap #3$}}

\newcommand{\relarr}
 {\looparrowright}

\newcommand{\rlim}
 {{\mbox{$\lower.95ex\hbox{{\rm lim}}$}\atop{\scriptstyle
{\rightarrow}}}{}}

\newcommand{\rlimd}
 {\lower0.37ex\hbox{$\pile{\lim \\ {\scriptstyle
\rightarrow}}$}{}}

\newcommand{\rstp}[3]
 {\mbox{$#1\! : #2 \relarr #3$}}

\newcommand{\scn}
 {{\bf scn}}

\newcommand{\scnS}
 {{\bf scn}_{\cal S}}

\newcommand{\semid}
 {\rtimes}

\newcommand{\Sep}
 {{\bf Sep}}

\newcommand{\sep}
 {{\bf sep}}

\newcommand{\Set}
 {{\bf Set }}

\newcommand{\Sh}
 {{\bf Sh}}

\newcommand{\ShE}
 {{\bf Sh}_{\cal E}}

\newcommand{\ShS}
 {{\bf Sh}_{\cal S}}

\newcommand{\sh}
 {{\bf sh}}

\newcommand{\Simp}
 {{\bf \Delta}}

\newcommand{\snd}
 {{\sf snd}}

\newcommand{\stg}[1]
 {\vec{#1}}

\newcommand{\stp}[3]
 {\mbox{$#1\! : #2 \to #3$}}

\newcommand{\Sub}
 {{\rm Sub}}

\newcommand{\SUB}
 {{\mathbb S}{\sf ub}}

\newcommand{\tbel}
 {\prec\!\prec}

\newcommand{\tic}[2]
 {\mbox{$#1\!.\!#2$}}

\newcommand{\tp}
 {\!:}

\newcommand{\tps}
 {:}

\newcommand{\tsub}
 {\pile{\lower0.5ex\hbox{.} \\ -}}

\newcommand{\wavy}
 {\leadsto}

\newcommand{\wavydown}
 {\,{\mbox{\raise.2ex\hbox{\hbox{$\wr$}
\kern-.73em{\lower.5ex\hbox{$\scriptstyle{\vee}$}}}}}\,}

\newcommand{\wbel}
 {\lt\!\lt}

\newcommand{\wstp}[3]
 {\mbox{$#1\!: #2 \wavy #3$}}

\newcommand{\fu}[2]
{[#1,#2]}

\newcommand{\st}[2]
 {\mbox{$#1 \to #2$}}

\section{The De Morgan topology}
In this section we first introduce, in the context of elementary toposes, the notion of De Morgan topology; this is shown to play, with respect to De Morgan toposes, the same role that the well-known notion of double-negation topology plays with respect to Boolean toposes. Then we provide explicit descriptions of the De Morgan topology on a presheaf topos (in terms of the corresponding Grothendieck topology) and on a localic topos (in terms of the corresponding locale).\\
Let us recall the following definitions.
\begin{definition}
An Heyting algebra $H$ is said to be a De Morgan algebra if and only if for each $p\in H$, $\neg p \vee \neg \neg p=1$.
\end{definition}
\begin{definition}
An elementary topos $\cal{E}$ is said to satisfy De Morgan's law (equivalently, to be a De Morgan topos) if its subobject classifier $\Omega$ is an internal De Morgan algebra in $\cal E$.
\end{definition}
There are many different known characterizations of De Morgan's law in toposes (we refer the reader to \cite{El2} for a comprehensive treatment); we will make use of the following ones:\\
(1) A topos $\cal E$ is De Morgan if and only if the canonical monomorphism $(\top, \bot):2=1\amalg 1\rightarrowtail \Omega_{\neg \neg}$ is an isomorphism in $\cal E$ (here $\Omega_{\neg \neg}$ denotes the equalizer in $\cal E$ of the pair of arrows $1_{\Omega}: \Omega \rightarrow \Omega$ and $\neg\neg: \Omega \rightarrow \Omega$, that is the subobject classifier of the topos $\sh_{\neg \neg}(\cal E)$).\\    
(2) A topos $\cal E$ is De Morgan if and only if the arrow $\bot:1\rightarrow \Omega$ (that is, the classifying map of the least subobject $0:0\rightarrowtail 1$ of $\Sub(1)$) has a complement in the Heyting algebra $\Sub(\Omega)$.\\
Let us recall from \cite{blasce} (or from \cite{El}) that, given an elementary topos $\cal E$ and a topology $j$ on it such that $j\leq \neg \neg$, $\sh_{j}(\cal E)$ is Boolean if and only if $j=\neg \neg$; in other words, $\neg \neg$ is the least topology $j$ on $\cal E$ among those which satisfy $j \leq \neg\neg$ such that $\sh_{j}(\cal E)$ is Boolean.\\

\begin{theorem}\label{thmDeMorgan}
Let $\cal E$ be an elementary topos. Then there exists a topology $m$ on $\cal E$ which is the least topology $j$ on $\cal E$ among those which satisfy $j \leq \neg\neg$ such that $\sh_{j}(\cal E)$ is De Morgan.\\  
\end{theorem}
The unique topology $m$ on $\cal E$ satisfying the condition of the proposition will be called the \emph{De Morgan topology} on $\cal E$, and the topos $\sh_{m}(\cal E)$ will be called the \emph{DeMorganization} of the topos $\cal E$.\\ 

\begin{proofs}
First, a remark on notation: given a local operator (i.e. a topology) $j$ on $\cal E$, we denote by $\Omega_{j}$ the equalizer in $\cal E$ of the pair of arrows $1_{\Omega}: \Omega \rightarrow \Omega$ and $j: \Omega \rightarrow \Omega$, that is the subobject classifier of the topos $\sh_{j}(\cal E)$, and by $a_{j}:{\cal E}\rightarrow \sh_{j}(\cal E)$ the associated sheaf functor.\\   
We define $m$ to be the smallest local operator $j$ on $\cal E$ such the canonical monomorphism $(\top, \bot):2=1\amalg 1\rightarrowtail \Omega_{\neg \neg}$ is $j$-dense; such an operator exists by a theorem of A. Joyal's (see Example A4.5.14 (b) p. 215 \cite{El}).\\
If $j$ is a topology such that $j \leq \neg\neg$ then $\sh_{\neg\neg}({\cal E})\subseteq \sh_{j}(\cal E)$ and $\Omega_{\neg\neg}\leq \Omega_{j}$ in $\Sub(\Omega)$; from this it easily follows that $(\Omega_{j})_{\neg\neg}=\Omega_{\neg\neg}$, equivalently $\sh_{\neg\neg}(\sh_{j}({\cal E}))=\sh_{\neg\neg}(\cal E)$ (cfr. the proof of Lemma A4.5.21 p. 220 \cite{El}). Since $\Omega_{\neg\neg}$ is a $j$-sheaf and the equality $(\Omega_{j})_{\neg\neg}=\Omega_{\neg\neg}$ holds, the fact the associated sheaf functor $a_{j}:{\cal E}\rightarrow \sh_{j}(\cal E)$ preserves coproducts implies that the canonical monomorphism $(\top, \bot):2=1\amalg 1\rightarrowtail {(\Omega_{j})}_{\neg \neg}$ for the topos $\sh_{j}(\cal E)$ can be obtained as the result of applying $a_{j}$ to the canonical monomorphism $(\top, \bot):2=1\amalg 1\rightarrowtail \Omega_{\neg \neg}$ for the topos $\cal E$. Hence, recalling that a monomorphism $m$ is $j$-dense for a local operator $j$ on $\cal E$ if and only if $a_{j}(m)$ is an isomorphism in $\sh_{j}(\cal E)$, we conclude that $m$ satisfies the following property: for each local operator $j$ such that $j\leq \neg\neg$,  $\sh_{j}(\cal E)$ is De Morgan if and only if $m\leq j$; this in particular implies that $m\leq \neg\neg$ (as $\sh_{\neg\neg}(\cal E)$ is always De Morgan being Boolean) and hence our thesis.  
\end{proofs}
Let us prove an analogous characterization for the double-negation topology.
\begin{theorem}
Let $\cal E$ be an elementary topos. Then the double-negation topology $\neg\neg$ on $\cal E$ is the least topology $j$ on $\cal E$ such that the canonical monomorphism $(\top, \bot):2=1\amalg 1\rightarrowtail \Omega$ is $j$-dense.\\
\end{theorem}
\begin{proofs}
Let us denote by $b$ the smallest local operator $j$ on $\cal E$ such that $(\top, \bot):2=1\amalg 1\rightarrowtail \Omega$ is $j$-dense; again, such an operator exists by Example A4.5.14 (b) p. 215 \cite{El}.\\ 
The canonical morphism $a_{\neg\neg}(\Omega)\rightarrow \Omega_{\neg\neg}$ is an isomorphism; indeed, this follows from Proposition A4.5.8 \cite{El} in view of the fact that the identity $\neg\neg(\neg\neg h \vee h)=1$ holds in any Heyting algebra. Hence, since $a_{\neg\neg}$ preserves coproducts, $a_{\neg\neg}((\top, \bot))$ is an isomorphism if and only if the canonical monomorphism $(\top, \bot)$ for the topos $\sh_{\neg\neg}({\cal E})$ is an isomorphism, and this is the case since $\sh_{\neg\neg}({\cal E})$ is Boolean. So we have that $(\top, \bot)$ is $\neg\neg$-dense and hence $b\leq \neg\neg$. Now, if $j\leq \neg\neg$ then $(\top, \bot)$ factors through $\Omega_{j}\mono \Omega$, so if $(\top, \bot)$ is $j$-dense then the factorization $(\top, \bot):2=1\amalg 1\rightarrowtail \Omega_{j}$ is $j$-dense (recall that the composite of two monomorphisms is dense with respect to a topology if and only if both of them are), in other words, $\sh_{j}({\cal E})$ is Boolean; in particular, $\sh_{b}({\cal E})$ is Boolean (as we have observed above that $b\leq \neg\neg$). Now, the facts that $b\leq \neg\neg$ and $\sh_{b}({\cal E})$ is Boolean together imply that $b=\neg\neg$, by the remark before Theorem \ref{thmDeMorgan}.\\
Note that it is possible to avoid invoking the existence of the topology $b$ in this proof by arguing as follows. If $j\leq \neg\neg$ then (we have observed above that) $(\top, \bot)$ $j$-dense implies $\sh_{j}(\cal E)$ Boolean, that is $j=\neg\neg$. For a general $j$, consider the meet $j\wedge \neg\neg$ in the lattice of topologies on $\cal E$. From the fact that meets in this lattice are computed ``pointwise'' and $(\top, \bot)$ is $\neg\neg$-dense (which we have observed above), we have that if $(\top, \bot)$ is $j$-dense then $(\top, \bot)$ is $(j\wedge \neg\neg)$-dense; so, since $j\wedge \neg\neg\leq \neg\neg$, we can refer to the previous case and conclude that $j\wedge \neg\neg=\neg\neg$ (equivalently, $\neg\neg \leq j$).
\end{proofs}

The following proposition states a couple of useful facts on the De Morgan topology.
\begin{proposition}\label{propDeMorgan}
Let $\cal E$ be an elementary topos and $m$ the De Morgan topology on it. Then\\
(i) $\sh_{m}({\cal E})={\cal E}$ if and only if $\cal E$ is a De Morgan topos;\\
(ii) For any topology $j$ on $\cal E$ such that $j\leq \neg\neg$, $\sh_{m}(\sh_{j}({\cal E}))=\sh_{k}(\cal E)$ where $k=m\vee j$ in the lattice of topologies on $\cal E$. 
\end{proposition}
\begin{proofs}
Part (i) is an immediate consequence of Theorem \ref{thmDeMorgan} as ${\cal E}=\sh_{1}(\cal E)$, where $1$ is the smallest topology on $\cal E$. Let us then prove part (ii). 
To prove our equality, we verify that the topology $k=m\vee j$ on $\cal E$ satisfies the universal property of the De Morgan topology on $\sh_{j}({\cal E})$ given by Theorem \ref{thmDeMorgan}. For a given topology $l$ on $\sh_{j}({\cal E})$ such that $l\leq\neg\neg_{\sh_{j}({\cal E})}$, $\sh_{l}(\sh_{j}({\cal E}))$ is a dense subtopos of $\sh_{j}({\cal E})$; but $\sh_{j}({\cal E})$ is a dense subtopos of $\cal E$ by hypothesis, so $\sh_{l}(\sh_{j}({\cal E}))$ is a dense subtopos of $\cal E$ (as the composition of dense inclusions is again a dense inclusion); then, by definition of De Morgan topology on $\cal E$, we have that $\sh_{l}(\sh_{j}({\cal E}))$ is De Morgan if and only if $\sh_{l}(\sh_{j}({\cal E}))\subseteq \sh_{m}(
{\cal E})$, if and only if $\sh_{l}(\sh_{j}({\cal E}))\subseteq \sh_{m}({\cal E})\cap \sh_{j}({\cal E})=\sh_{m\vee j}(\cal E)$.
\end{proofs}

Now, our aim is to describe explicitly the De Morgan topology on a presheaf topos $[{\cal C}^{\textrm{op}}, \Set]$. To this end, we rephrase criterion (2) above for a topos to be De Morgan in the case of a Grothendieck topos ${\cal E}=\Sh({\cal C}, J)$.\\
Recall that the subobject classifier $\Omega_{J}:{\cal C}^{\textrm{op}}\rightarrow \Set$ of the topos $\Sh({\cal C}, J)$ is defined by:\\
$\Omega_{J}(c)=\textrm{\{$R$ | $R$ is a $J$-closed sieve on $c$\}}$ (for an object $c\in \cal C$),\\
$\Omega_{J}(f)=f^{\ast}(-)$ (for an arrow $f$ in $\cal C$),\\
where $f^{\ast}(-)$ denotes the operation of pullback of sieves in $\cal C$ along $f$.\\ The arrow $\bot:1\rightarrow \Omega_{J}$ is the classifying map of the smallest subobject $0:0\rightarrow 1$ in $\Sub_{\cal E}(1)$, which is the subfunctor of $1$ defined by: $0(c)=1(c)=\{\ast\}$ if $\emptyset\in J(c)$ and $0(c)=\emptyset$ if $\emptyset \notin J(c)$. A formula p. 142 \cite{MM} then gives:
\[
\bot(c)(\ast)=\{f:d\rightarrow c \textrm{ | } \ast \in 0(d)\}=\{f:d\rightarrow c \textrm{ | } \emptyset \in J(d)\}.
\] 
(despite the notation, here and below the domains of the arrows are intended to be variable).\\
Let us put for convenience $R_{c}:=\{f:d\rightarrow c \textrm{ | } \emptyset \in J(d)\}$, for $c\in \cal C$.\\
By using formula (19) p. 149 \cite{MM} we get
\[
(\neg \bot)(c)=\{R\in \Omega_{J}(c) \textrm{ | for any } f:d\rightarrow c, f^{\ast}(R)=R_{d}\textrm{ implies } f\in R_{c} \},
\]
for any $c\in {\cal C}$.\\
Let us now calculate $\bot\vee \neg \bot$ by using formula (5) p. 145 \cite{MM}:
\[
(\bot\vee \neg \bot)(c)=\{R\in \Omega_{J}(c) \textrm{ | } \{f:d\rightarrow c \textrm{ | } f^{\ast}(R)=R_{d} \textrm{ or } f^{\ast}(R)\in (\neg \bot)(d)\}\in J(c)\},
\]
for any $c\in {\cal C}$.\\
Hence we conclude that $\bot\vee \neg \bot=1_{\Omega_{J}}$ (equivalently, $\Sh({\cal C}, J)$ is De Morgan) if and only if for every object $c\in \cal C$ and $J$-closed sieve $R$ on $c$
\[
\{f:d\rightarrow c \textrm{ | } (f^{\ast}(R)=R_{d}) \textrm{ or } (\textrm{for any } g:e\rightarrow d, g^{\ast}(f^{\ast}(R))=R_{e}\textrm{ implies } g\in R_{d})\}
\]
belongs to $J(c)$.\\
Let us now restrict our attention to Grothendieck topologies $J$ on $\cal C$ such that all $J$-covering sieves are non-empty. Under this hypothesis, we have that $R_{c}=\emptyset$ (for each $c\in \cal C$) and hence
\[
(\neg \bot)(c)=\{R\in \Omega_{J}(c) \textrm{ | for any } f:d\rightarrow c, f^{\ast}(R)\neq \emptyset\},
\]
This motivates the following definition: a sieve $R$ on $c\in \cal C$ is said to be stably non-empty if for any $f:d\rightarrow c, f^{\ast}(R)\neq \emptyset$.\\
Let us put, for any sieve $R$ on $c\in \cal C$,
\[
M_{R}:=\{f:d\rightarrow c \textrm{ | } (f^{\ast}(R)=\emptyset) \textrm{ or } (f^{\ast}(R) \textrm{ is stably non-empty}) \}.
\]
Then we have
\[
(\bot\vee \neg \bot)(c)=\{R\in \Omega_{J}(c) \textrm{ | } M_{R}\in J(c)\}.
\]
So, under the hypothesis that every $J$-covering sieve in non-empty, we get the following simplified form of our criterion:
$\Sh({\cal C}, J)$ is De Morgan if and only if for every object $c\in \cal C$ and $J$-closed sieve $R$ on $c$, $M_{R}\in J(c)$.
\begin{remarks}\label{rmks}

(a) For any sieve $R$ on $c\in \cal C$ and any arrow $f:d\rightarrow c$ in $\cal C$, $f^{\ast}(M_{R})=M_{f^{\ast}(R)}$.\\
(b) If $r:c'\rightarrow c$ is a monomorphism in $\cal C$ then, given a sieve $R'$ on $c'$ and denoted by $R$ the sieve $\{r\circ f \textrm{ | } f\in R'\}$ on $c$, we have that $r^{\ast}(M_{R})=M_{R'}$. Indeed, by (a) we have $r^{\ast}(M_{R})=M_{r^{\ast}(R)}$, so it is enough to prove that $r^{\ast}(R)=R'$; one inclusion is obvious, while the other holds since $r'$ is monic. This implies that if $J$ is a Grothendieck topology on $\cal C$ then $M_{R}\in J(c)$ implies $M_{R'}\in J(c')$; thus, in checking that the condition of our criterion is satisfied, we can restrict our attention to any collection $\cal F$ of objects in $\cal C$ with the property that for each object $c$ in ${\cal C}$ there is a monomorphism $r$ in ${\cal C}$ from $c$ to an object in $\cal F$.\\
(c) Under the hypothesis that every $J$-covering sieve in non-empty, if $\overline{R}$ is the $J$-closure of the sieve $R$, $M_{\overline{R}}=M_{R}$; indeed, a sieve $R$ is empty if and only if its $J$-closure $\overline{R}$ is.   
\end{remarks}       

Remark \ref{rmks}(c) implies that in the simplified form of our criterion above we can equivalently quantify over \emph{all} sieves $R$ in $\cal C$. This leads us to give the following definition: given a category $\cal C$, the De Morgan topology $M_{\cal C}$ on it is the Grothendieck topology on $\cal C$ generated by the family of sieves $\{M_{R} \textrm{ | $R$ sieve in } {\cal C}\}$. In fact, our criterion says precisely that, for any Grothendieck topology $J$ such that every $J$-covering sieve is non-empty (equivalently, $J\leq \neg\neg_{[{\cal C}^{\textrm{op}},\Set]}$), $\Sh({\cal C}, J)$ is De Morgan if and only $M_{\cal C}\leq J$. This, together with the observation that every $M_{\cal C}$-covering sieve is non-empty, proves that the De Morgan topology $M_{\cal C}$ on $\cal C$ is exactly the Grothendieck topology on $\cal C$ corresponding to the De Morgan topology on the topos $[{\cal C}^{\textrm{op}},\Set]$.\\
Summarizing, we have the following result. 
\begin{theorem}\label{criterion}
Let $\cal C$ be a category. There exists a Grothendieck topology $M_{\cal C}$ on $\cal C$, called the De Morgan topology on $\cal C$, which satisfies the following property: for any Grothendieck topology $J$ such that every $J$-covering sieve is non-empty, $\Sh({\cal C}, J)$ is De Morgan if and only $M_{\cal C}\leq J$. $M_{\cal C}$ is the Grothendieck topology on $\cal C$ corresponding to the De Morgan topology on the topos $[{\cal C}^{\textrm{op}},\Set]$ and is generated by the family of sieves $\{M_{R} \textrm{ $|$ R sieve in } {\cal C}\}$.
\end{theorem}\qed 
Notice that in case $J$ is the trivial topology the theorem immediately gives the following characterization: the topos $[{\cal C}^{\textrm{op}},\Set]$ is De Morgan if and only if $\cal C$ satisfies the right Ore condition (this is a well-known result, cfr. Example D4.6.3(a) p. 1001 \cite{El2}). Also, if $\cal C$ satisfies the right Ore condition, then $M_{\cal C}$ is clearly the trivial Grothendieck topology on $\cal C$; as a consequence, we obtain the following result.
\begin{corollary}\label{corfond}
Let $\cal C$ be a category satisfying the right Ore condition. Then for every Grothendieck topology $J$ such that every $J$-covering sieve is non-empty, $\Sh({\cal C}, J)$ is De Morgan.  
\end{corollary}\qed 

Now, let us briefly turn our attention to Boolean toposes. Starting from the well-known characterization: $\Sh({\cal C}, J)$ is Boolean if and only if $\bot \vee \top=1_{\Omega_{J}}$ in $\Sub(\Omega_{J})$, our methods can be easily adapted to prove
the following criterion:\\
$\Sh({\cal C}, J)$ is Boolean if and only if for every object $c\in \cal C$ and $J$-closed sieve $R$ on $c$,
\[
\{f:d\rightarrow c \textrm{ | } (f^{\ast}(R)=R_{d}) \textrm{ or } (f\in R)\}\in J(c). 
\]   
If $J$ is a Grothendieck topology on $\cal C$ such that every $J$-covering sieve is non-empty, the criterion becomes:\\$\Sh({\cal C}, J)$ is Boolean if and only if for every object $c\in \cal C$ and $J$-closed sieve $R$ on $c$, $B_{R}:=\{f:d\rightarrow c \textrm{ | } (f^{\ast}(R)=\emptyset) \textrm{ or } (f\in R)\}\in J(c)$.\\ In fact, the condition `$J$-closed' here can be put in parentheses, by the following characherization of the double-negation topology $\neg\neg_{[{\cal C}^{\textrm{op}},\Set]}$ on the topos $[{\cal C}^{\textrm{op}},\Set]$ (and the remark preceding Theorem \ref{thmDeMorgan}).\\ 
Recall from \cite{MM} that $\neg\neg_{[{\cal C}^{\textrm{op}},\Set]}$ corresponds to the dense topology $D$ on $\cal C$, that is to the Grothendieck topology $D$ on $\cal C$ defined as follows: for a sieve $R$ in $\cal C$
\[
R\in D(c)\textrm{ iff $R$ is stably non-empty}.  
\]
It is immediate to prove that the topology $D$ is given precisely by the collection of sieves $\{B_{R} \textrm{ | $R$ sieve in } {\cal C}\}$; indeed, for any sieve $R$ on $c\in \cal C$, $B_{R}\in D(c)$ and for any $R\in D(c), R=B_{R}$.\\ 
Suppose now that $\cal C$ satisfies the right Ore condition and every $J$-covering sieve is non-empty; the criterion above simplifies to:\\
\[
\begin{array}{ccl}
\Sh({\cal C}, J) \textrm{ is Boolean} & \textrm{iff} & \textrm{for every object $c\in \cal C$ and $J$-closed sieve $R$ on $c$,}\\
& & R\cup \{f:d\rightarrow c \textrm{ | } f^{\ast}(R)=\emptyset\}\in J(c),\\
& \textrm{iff} & \textrm{every non-empty $J$-closed sieve is $J$-covering,}\\
& \textrm{iff} & \textrm{the only non-empty $J$-closed sieves are the maximal sieves.}   
\end{array}
\]
Finally, let us describe the De Morgan topology on a given topos $\Sh(X)$ of sheaves a locale $X$.
To this end, we prove the following result, which is the natural embodiment of a number of ideas present in \cite{El} and \cite{El2} (the notation used below being that of \cite{El} and \cite{El2}). 
\begin{proposition}\label{loctop}
Let $X$ be a locale. Then there exists a frame isomorphism $N({\cal O}(X))\cong {\bf Lop}(\Sh(X))$ from the frame $N({\cal O}(X))$ of nuclei on the frame ${\cal O}(X)$ corresponding to $X$ and the frame ${\bf Lop}(\Sh(X))$ of local operators on the topos $\Sh(X)$ (equivalently, a coframe isomorphism between the coframe of sublocales of $X$ and the coframe of subtoposes of $\Sh(X)$). Through this isomorphism, an open (resp. closed) nucleus on an element $a\in {\cal O}(X)$ corresponds to the open (resp. closed) subtopos of $\Sh(X)$ determined by $a$ (regarded as a subterminal object in \Sh(X)), and the dense-closed factorization of a given sublocale of $X$ corresponds to the dense-closed factorization of the corresponding geometric inclusion.   
\end{proposition} 
\begin{proofs}
Given a geometric inclusion $i:{\cal E}\rightarrow \Sh(X)$ with codomain $\Sh(X)$, $i$ is localic (cfr. Example A4.6.2(a) \cite{El}), hence the topos $\cal E$ is localic (cfr. Example A4.6.2(e) \cite{El} and Theorem C1.4.7 \cite{El2}), that is there exists a locale $Y$ such that ${\cal E}\simeq \Sh(Y)$; by Proposition C1.4.5 \cite{El2} such a locale $Y$ is unique up to isomorphism in the category ${\bf Loc}$ of locales and by Corollary C1.5.2 \cite{El2} the inclusion $i:\Sh(Y)\rightarrow \Sh(X)$ corresponds to a unique subobject $Y\rightarrow X$ in ${\bf Loc}$. Conversely, any sublocale $L$ of $X$ gives rise to a geometric inclusion $\Sh(L)\rightarrow \Sh(X)$ (again by Corollary C1.5.2 \cite{El2}). These two assignments are clearly inverse to each other, and hence define a bijection between the (equivalence classes of) geometric inclusions with codomain $\Sh(X)$ and the sublocales of $X$. Now, recalling that the (equivalence classes of) geometric inclusions with codomain $\Sh(X)$ are in bijection with the local operators on the topos $\Sh(X)$ and the sublocales of $X$ are in bijection with the nuclei on the frame ${\cal O}(X)$, we obtain a bijection $N({\cal O}(X))\cong {\bf Lop}(\Sh(X))$. This bijection is in fact a frame isomorphism; indeed, given two sublocales of $X$, they are included one into the other if and only if the corresponding subtoposes are (again, this is an immediate consequence of Proposition C1.4.5 \cite{El2} and Corollary C1.5.2 \cite{El2}). This concludes the proof of the first part of the proposition. Now, if $o(a)$ is the open nucleus on an element $a\in {\cal O}(X)$ then the subtopos corresponding to it via the isomorphism is the open subtopos $\tilde{o}(a)$ determined by $a$, $a$ being regarded here as a subterminal object in \Sh(X) (cfr. the discussion p. 204 \cite{El}); from this we deduce that the closed nucleus $c(a)$ on an element $a\in {\cal O}(X)$ corresponds to the closed subtopos $\tilde{c}(a)$ determined by $a$, as $c(a)$ and $\tilde{c}(a)$ are respectively the complements of $o(a)$ and $\tilde{o}(a)$ in the frames $N({\cal O}(X))$ and ${\bf Lop}(\Sh(X))$ (cfr. section A4.5 \cite{El} and Example C1.1.16(b) \cite{El2}).\\
Recall from \cite{El2} that every sublocale $Y$ of a given locale $X$ has a closure $\overline{Y}$; specifically, if $j$ is the nucleus on ${\cal O}(X)$ corresponding to $Y$ then $c(j(0))$ is the nucleus on ${\cal O}(X)$ corresponding to $\overline{Y}$. In passing, we note that $\overline{Y}$ is characterized among the sublocales of $X$ by the following property: it is the largest sublocale $Z$ of $X$ such that for each open sublocale $A$ of $X$ $A\cap Z\neq \emptyset$ (if) and only if $A\cap Y\neq \emptyset$; indeed, by considering the corresponding fixsets, it is immediate to see that $A\cap Y\neq \emptyset$ if and only if $a\geq j(0)$, where $A=o(a)$. Then we have a factorization $Y\rightarrow \overline{Y}\rightarrow X$ where $Y\rightarrow \overline{Y}$ is dense and $\overline{Y}\rightarrow X$ is closed. We want to show that the corresponing geometric inclusions $\Sh(Y)\rightarrow \Sh(\overline{Y})$ and $\Sh(\overline{Y}) \rightarrow \Sh(X)$ are respectively dense and closed. We recall from \cite{El} that the dense-closed factorization of a geometric inclusion $\sh_{j}(\cal E)\hookrightarrow {\cal E}$ is given by $\sh_{j}({\cal E})\rightarrow \sh_{\tilde{c}(ext(j))}({\cal E})\rightarrow {\cal E}$, where $\textrm{ext}: {\bf Lop}({\cal E})\rightarrow \Sub_{\cal E}(1)$ is the right adjoint to the map $\tilde{c}:\Sub_{\cal E}(1)\rightarrow {\bf Lop}(\cal E)$ sending each subterminal object to the corresponding closed subtopos. Now, if ${\cal E}=\Sh(X)$ the map $\tilde{c}$ corresponds via our isomorphism to the map $c:{\cal O}(X)\rightarrow N({\cal O}(X))$ sending an element $a\in {\cal O}(X)$ to the corresponding closed nucleus $c(a)$; by arguing in terms of fixsets, it is immediate to verify that this map has a right adjoint given by the map sending a nucleus $j$ to its value $j(0)$ at 0; hence, the dense-closed factorization of the subtopos $\Sh(Y) \rightarrow \Sh(X)$ is given by $\Sh(Y)\rightarrow \Sh(\overline{Y}) \rightarrow \Sh(X)$, as required.       \end{proofs}
We are now ready to solve our original problem.
\begin{theorem}
Let $X$ be a locale. Then the DeMorganization $\Sh_{m}(\Sh(X))$ of the topos $\Sh(X)$ is equivalent to the topos $\Sh(X_{m})$ of sheaves on the locale $X_{m}$ defined as follows: ${\cal O}(X_{m})$ is the quotient of ${\cal O}(X)$ by the filter generated by the family $\{u\vee\neg u \textrm{ | $u$ is a regular element of } {\cal O}(X)\}$.
\end{theorem}
\begin{proofs}
By definition of De Morgan topology, $\Sh_{m}(\Sh(X))$ is the largest dense De Morgan subtopos of $\Sh(X)$. In view of Proposition \ref{loctop} and of the well-known characterization `$\Sh(X)$ is a De morgan topos if and only if $X$ is a De Morgan locale (i.e. ${\cal O}(X)$ is a De Morgan algebra)', it is equivalent to prove that $X_{m}$ is the largest dense De Morgan sublocale of $X$. This will immediately follow from the definition of $X_{m}$, once we have proved that $X_{m}$ is dense in $X$. Indeed, if $L$ is a sublocale of $X$ with corresponding surjective homomorphism of frames $l:{\cal O}(X)\rightarrow {\cal O}(L)$ then $L$ is dense in $X$ if and only if for each $a\in {\cal O}(X)$, $l(a)=0$ implies $a=0$; so, if $L$ is dense in $X$, $l$ preserves the operation of pseudocomplementation and hence $L$ is a De Morgan locale if and only if $l$ factors through the natural projection ${\cal O}(X)\rightarrow {\cal O}(X_{m})$.\\
Now, if $j$ is the nucleus corresponding to the sublocale $X_{m}$ then to prove that $X_{m}$ is dense amounts to verify that $j(0)=0$. By definition of nucleus corresponding to (a sublocale regarded as) a surjective homomorphism of frames, we have that $j(0)$ is the largest element $a\in {\cal O}(X)$ such that both $a\imp 0$ and $0\imp a$ belongs to the filter in the statement of the proposition; so we have to prove that for any $a\in {\cal O}(X)$, $\neg a$ belongs to the filter if and only if $a=0$. For $a$ to belong to the filter it is necessary (and sufficient) that there exists a finite number $u_{1}, u_{2}, \ldots, u_{n}$ of regular elements of ${\cal O}(X)$ such that $a\geq \mathbin{\mathop{\textrm{\huge $\vee$}}\limits_{1\leq i\leq n}}(u_{i}\vee \neg u_{i})$. Now, denoted by ${\cal O}(X)_{\neg\neg}$ the lattice of regular elements of ${\cal O}(X)$, the double negation operator $\neg\neg$ is a frame homomorphism ${\cal O}(X)\rightarrow {\cal O}(X)_{\neg\neg}$ and hence by applying it to the inequality above we obtain that $\neg a=\neg\neg\neg a$ is the top element of the Boolean algebra ${\cal O}(X)_{\neg\neg}$, equivalently $a=0$.          
\end{proofs}

\section{The simplification method} 
The purpose of this section is to give a simplified description of our criterion for a Grothendieck topos to be De Morgan, and in particular of the De Morgan topology, in several cases of interest.\\
Let us start with an informal description of our strategy. The main idea is that the more categorical stucture we have on $\cal C$, the more we should be able to simplify the description of our criterion. This simplification will in fact be carried out in three steps; at each step the category $\cal C$ will be supposed to have some more categorical structure than it had in the previous step and, as a result, a simpler description of the criterion will be achieved.\\
Let $({\cal C}, J)$ be a Grothendieck site. Then, denoted by $\tilde{{\cal C}}$ the full subcategory of $\cal C$ on the objects which are not $J$-covered by the empty sieve and by $\tilde{J}$ the topology induced by $J$ on $\tilde{{\cal C}}$, the toposes $\Sh({\cal C}, J)$ and $\Sh(\tilde{{\cal C}}, \tilde{J})$ are naturally equivalent (cfr. Example C2.2.4(e) \cite{El2}). Theorem \ref{criterion} then implies that a Grothendieck topos $\Sh({\cal C}, J)$ is De Morgan if and only if $M_{\tilde{{\cal C}}}\leq \tilde{J}$. In investigating whether a Grothendieck topos $\Sh({\cal C}, J)$ is De Morgan, we would then naturally opt for using, because of its simplicity, this latter form of the criterion which involves working with the category $\tilde{\cal C}$ rather than with $\cal C$. However, while our original category ${\cal C}$ may have a certain amount of categorical structure, by passing from $\cal C$ to $\tilde{\cal C}$ it often happens that a lot of categorical structure is lost. Our stategy will be then to work with the site $(\tilde{{\cal C}}, \tilde{J})$, but by keeping in mind its relationship with the original site $({\cal C}, J)$ (what we exactly mean by this will be clear later). We will restrict our attention to Grothendieck topologies $J$ such that the only object of $\cal C$ which is $J$-covered by the empty sieve is the initial object $0_{\cal C}$ (up to isomorphism) (note that for a subcanonical topology $J$, this is always the case). Also, instead of requiring that the category $\cal C$ has enough structure itself, we will more loosely require $\cal C$ to be closed (in the obvious sense) under the categorical structure on a larger category $\cal D$; that is, we will work in the context of (full) embeddings $\tilde{\cal C}\hookrightarrow {\cal C} \hookrightarrow {\cal D}$, where $\cal D$ is supposed to be a category ``with enough structure'' and $\cal C$ is assumed to be closed under this structure.\\ 
First, let us introduce some terminology.\\   
Given an embedding ${\cal C}\hookrightarrow {\cal D}$, where $\cal D$ is a category with pullbacks, and two arrows $f:a\rightarrow c$ and $g:b\rightarrow c$ in $\cal D$ with common codomain, we denote by $\textrm{p.b.($f$, $g$)}$ the object $p$ in $\cal D$ forming the pullback square
\[  
\xymatrix {
p \ar[d] \ar[r] & a \ar[d]^{f}  \\
b \ar[r]_{g} & c }
\]\\  
in $\cal D$; of course, $p$ is defined only up to isomorphism in $\cal D$.\\
The following proposition represents the first step of our simplification process. Below, for $\cal C$ to be closed in $\cal D$ under pullbacks we mean that whenever we have a pullback square
\[  
\xymatrix {
p \ar[d] \ar[r] & a \ar[d]^{f}  \\
b \ar[r]_{g} & c }
\]\\  
in $\cal D$ where $f$ and $g$ lie in $\cal C$ then (an isomorphic copy of) the object $p$, and hence the whole square, also lies in $\cal C$. 
 
\begin{proposition}\label{firststep}
Let ${\cal C}\hookrightarrow {\cal D}$ be a full embedding of categories such that $\cal D$ has pullbacks and a strict initial object $0\in {\cal C}$ and $\cal C$ is closed in $\cal D$ under pullbacks. Then for any object $c\in \tilde{{\cal C}}$, sieve $R$ on $c$ in $\tilde{{\cal C}}$ and arrow $f:d\rightarrow c$ in $\tilde{{\cal C}}$ we have: 
\[
\begin{array}{ccl}
f^{\ast}(R)=\emptyset & \textrm{iff} & \textrm{for every arrow $r$ in $R$ }, \textrm{p.b.($f$, $r$)}\cong 0;\\
f^{\ast}(R) \textrm{ is stably non-empty}  & \textrm{iff} & \textrm{for every arrow $g$ in $\tilde{{\cal C}}$ s.t. $cod(g)=dom(f)$},\\
& & \textrm{there exists $r$ in $R$ with } \textrm{p.b.($g$, $r$)}\ncong 0.\\ 
\end{array}
\]
(the sieve pullbacks $f^{\ast}(R)$ above being taken in the category $\tilde{\cal C}$).
\end{proposition} 
\begin{proofs}
Let us prove the first assertion, the second being an immediate consequence of it.\\
Let us suppose that $f^{\ast}(R)=\emptyset$. If for an arrow $r$ in $R$ we had $\textrm{p.b.($f$, $r$)}\ncong 0$ then we would have a pullback square
\[  
\xymatrix {
p \ar[d]_{h} \ar[r]^{k} & d \ar[d]^{f}  \\
b \ar[r]_{r} & c }
\]\\  
in $\cal D$ with $P\ncong 0$; hence, since $\cal C$ is closed in $\cal D$ under pullbacks, the arrow $k$ would lie in $\tilde{\cal C}$ and satisfy $k\in f^{\ast}(R)$, contradicting our assumption.\\
Conversely, let us suppose that $f^{\ast}(R)$ is non-empty. Then there exists an arrow $k:e\rightarrow d$ in $\tilde{\cal C}$ such that $f\circ k$ belongs to $R$. Hence $e\ncong 0$ and we have a commutative square  
\[  
\xymatrix {
e \ar[d]_{1_{d}} \ar[r]^{k} & d \ar[d]^{f}  \\
e \ar[r]_{f\circ k} & c }
\]\\  
Then, $0$ being a strict initial object in $\cal D$, by the universal property of the pullback it follows that $\textrm{p.b.($f$, $f\circ k$)}\ncong 0$.       
\end{proofs}
It is sensible at this point to introduce the following terminology: given two arrows $f$ and $g$ in $\tilde{\cal C}$ with common codomain, they are said to be \emph{disjoint} (equivalently, $f$ is said to be disjoint from $g$) if $\textrm{p.b.($f$, $g$)}\cong 0$, while $f$ is said to be \emph{stably joint} with $g$ if for each arrow $k$ in $\tilde{\cal C}$ such that $cod(k)=dom(f)$ we have $\textrm{p.b.($f\circ k$, $g$)}\ncong 0$.\\
We note that Proposition \ref{firststep} implies that, given two arrows $f$ and $r$ in $\tilde{\cal C}$ with common codomain, $f^{\ast}((r))=\emptyset$ if and only if $f$ and $r$ are disjoint, while $f^{\ast}((r))$ is stably non-empty if and only if $f$ is stably joint with $r$.\\   

Let us go on to the second step. Below, for $\cal C$ to be closed in $\cal D$ under cover-mono factorizations we mean that if $d\epi c'\mono c$ is the cover-mono factorization in $\cal D$ of a morphism $d\rightarrow c$ lying in $\cal C$, then (an isomorphic copy of) the object $c'$ (and hence the whole factorization) also lies in $\cal C$. 
\begin{proposition}\label{secondstep}
Let ${\cal C}\hookrightarrow {\cal D}$ be a full embedding of categories such that $\cal D$ is a regular category having a strict initial object $0\in \cal C$ and $\cal C$ is closed in $\cal D$ under pullbacks and cover-mono factorizations. Given an object $c\in \tilde{{\cal C}}$, a sieve $R$ on $c$ in $\tilde{{\cal C}}$ and an arrow $f:d\rightarrow c$ in $\tilde{{\cal C}}$, let us denote, for each arrow $r$ in $R$, by $dom(r) \epi x \stackrel{r'}{\mono} c$ its cover-mono factorization in $\cal D$ and by $R'$ the sieve in $\tilde{\cal C}$ generated by the arrows $r'$ (for $r$ in $R$). Then we have: 
\[
\begin{array}{ccl}
f^{\ast}(R)=\emptyset & \textrm{iff} &  f^{\ast}(R')=\emptyset;\\
f^{\ast}(R) \textrm{ is stably non-empty}  & \textrm{iff} & f^{\ast}(R') \textrm{ is stably non-empty}.\\  
\end{array}
\]
\end{proposition}
\begin{proofs}
Of course, it is enough to prove the first equivalence. This easily follows from Proposition \ref{firststep} and our hypotheses. Indeed, we have that $f^{\ast}(R)=\emptyset$ if and only if for each $r$ in $R$ $\textrm{p.b.($f$, $r$)}\cong 0$, if and only if for each $r'$ in $R'$ $\textrm{p.b.($f$, $r$)}\cong 0$, if and only if $f^{\ast}(R')=\emptyset$, where the second equivalence follows from the fact that, given a cover $d\epi c$, $c\cong 0$ if and only if $d\cong 0$ ($0$ being a strict initial object). 
\end{proofs} 
\begin{corollary}\label{cor1}
Let $({\cal C}, J)$ be a Grothendieck site such that the only object of $\cal C$ which is $J$-covered by the empty sieve is the initial object $0$ of $\cal C$ (up to isomorphism) and ${\cal C}\hookrightarrow {\cal D}$ a full embedding of categories such that $\cal D$ is a regular category having a strict initial object $0\in \cal C$ and $\cal C$ is closed in $\cal D$ under pullbacks and cover-mono factorizations. Then $\Sh({\cal C}, J)$ is a De Morgan topos if and only if for each sieve $R$ in $\tilde{\cal C}$ generated in $\tilde{\cal C}$ by morphisms which are monic in $\cal D$, $M_{R}=\{f:d\rightarrow c \textrm{ in $\tilde{\cal C}$ | } (f^{\ast}(R)=\emptyset) \textrm{ or } (f^{\ast}(R) \textrm{ is stably non-empty}) \}$ is a $\tilde{J}$-covering sieve. 
\end{corollary}
\begin{proofs}
From Proposition \ref{secondstep} we have that $M_{R}=M_{R'}$; our thesis then follows from the remarks at the beginning of this section. 
\end{proofs}
Let us now proceed to the third step. Below, for $\cal C$ to be closed in $\cal D$ under arbitrary (i.e. set-indexed) unions of subobjects we mean that whenever we have a set of arrows in $\cal C$ with common codomain $c\in {\cal C}$ which are monic in $\cal D$, the union of them in $\Sub_{\cal D}(c)$ also lies (up to isomorphism) in $\cal C$.
\begin{proposition}\label{thirdstep}
Let ${\cal C}\hookrightarrow {\cal D}$ be a full embedding of categories such that $\cal D$ is a geometric category with a (strict) initial object $0\in \cal C$ and $\cal C$ is closed in $\cal D$ under pullbacks, cover-mono factorizations and arbitrary unions of subobjects. Given an object $c\in \tilde{{\cal C}}$, a sieve $R$ on $c$ in $\tilde{{\cal C}}$ generated by arrows $\{r_{i}, i\in I\}$ which are monic in $\cal D$, and an arrow $f:d\rightarrow c$ in $\tilde{{\cal C}}$, let us denote by $r$ the union of the subobjects $r_{i}$ (for $i\in I$) in $\Sub_{\cal D}(c)$ and by (r) the sieve generated by $r$ in $\tilde{\cal C}$. Then we have: 
\[
\begin{array}{ccl}
f^{\ast}(R)=\emptyset & \textrm{iff} &  f^{\ast}((r))=\emptyset;\\
f^{\ast}(R) \textrm{ is stably non-empty}  & \textrm{iff} & f^{\ast}((r)) \textrm{ is stably non-empty}.\\  
\end{array}
\]
\end{proposition}

\begin{proofs}
This immediately follows from Proposition \ref{firststep} and the fact that unions of subobjects in $\cal D$ are stable under pullback; indeed, we have that $f^{\ast}(R)=\emptyset$ if and only if for each $r_{i}$ in $R$ $\textrm{p.b.($f$, $r_{i}$)}\cong 0$, if and only if $\mathbin{\mathop{\textrm{\huge $\cup$}}\limits_{i\in I}}\textrm{p.b.($f$, $r_{i}$)}\cong 0$, if and only if $\textrm{p.b.($f$, $r$)}\cong 0$, if and only if $f^{\ast}((r))=\emptyset$.  
\end{proofs}
From Propositions \ref{secondstep} and \ref{thirdstep} we immediately deduce the following corollary.
\begin{corollary}\label{cor2}
Let $({\cal C}, J)$ be a Grothendieck site such that the only object of $\cal C$ which is $J$-covered by the empty sieve is the initial object $0$ of $\cal C$ (up to isomorphism) and ${\cal C}\hookrightarrow {\cal D}$ be a full embedding of categories such that $\cal D$ is a geometric category with a (strict) initial object $0\in \cal C$ and $\cal C$ is closed in $\cal D$ under pullbacks, cover-mono factorizations and arbitrary unions of subobjects. Then $\Sh({\cal C}, J)$ is a De Morgan topos if and only if for each arrow $r$ in $\tilde{\cal C}$ which is monic in $\cal D$ $M_{(r)}=\{f:d\rightarrow c \textrm{ in $\tilde{\cal C}$ | ($f$ is disjoint from $r$) or ($f$ is stably joint with $r$)} \}$\\ is a $\tilde{J}$-covering sieve. 
\end{corollary}\qed
The following propositions are the analogues ``for the arrows $f$'' of Propositions \ref{secondstep} and \ref{thirdstep}.
\newpage
\begin{proposition}\label{arrowsf1}
Let ${\cal C}\hookrightarrow {\cal D}$ be a full embedding of categories such that $\cal D$ is a regular category having a strict initial object $0\in \cal C$ and $\cal C$ is closed in $\cal D$ under pullbacks and cover-mono factorizations. Given an object $c\in \tilde{{\cal C}}$, a sieve $R$ on $c$ in $\tilde{{\cal C}}$ and an arrow $f:d\rightarrow c$ in $\tilde{{\cal C}}$, if $d \stackrel{f''}{\epi} x \stackrel{f'}{\mono} c$ is the cover-mono factorization of $f$ in $\cal D$ then
\[
\begin{array}{ccl}
f^{\ast}(R)=\emptyset & \textrm{iff} &  f'^{\ast}(R)=\emptyset;\\
f^{\ast}(R) \textrm{ is stably non-empty}  & \textrm{iff} & f'^{\ast}(R) \textrm{ is stably non-empty}.\\  
\end{array}
\]
\end{proposition}
\begin{proofs}
Let us begin to prove the first part of the proposition. One implication is obvious as $f$ factors through $f'$; let us prove the other one. Suppose that $f'^{\ast}(R)\neq\emptyset$. Then there exists an arrow $k:dom(k)\rightarrow x$ in $\tilde{\cal C}$ such that $f'\circ k$ belongs to $R$. Now, since $f'$ is monic $\textrm{p.b.($f'\circ k$, $f'$)}=dom(k)\ncong 0$ so $\textrm{p.b.($f'\circ k$, $f$)}\ncong 0$ as $f''$ is a cover and $0$ is strictly initial. This implies that $f^{\ast}(R)\neq\emptyset$ by Proposition \ref{firststep}. This concludes the proof of the first part.\\
Let us now prove the second part. Again, one direction is trivial. To prove the other implication, let us suppose that $f^{\ast}(R)$ is stably non-empty. Given any arrow $g:e\rightarrow x$ in $\tilde{\cal C}$, we want to prove, according to Proposition \ref{firststep}, that there exists an arrow $r$ in $R$ such that $\textrm{p.b.($f'\circ g$, $r$)}\ncong 0$. To find such an arrow $r$, consider in $\cal D$ the pullback
\[  
\xymatrix {
y \ar[d]_{h} \ar[r]^{k} & d \ar[d]^{f''}  \\
e \ar[r]_{g} & x }
\]\\  
As $f''$ is a cover then $h$ is a cover, so $y\ncong 0$ and $k$ is an arrow in $\tilde{\cal C}$; then, $f^{\ast}(R)$ being stably non-empty, there exists an arrow $r$ in $R$ such that $\textrm{p.b.($f\circ k$, $r$)}\ncong 0$. From this it is immediate to see (by using that $h$ is a cover and $0$ is strictly initial) that $\textrm{p.b.($f'\circ g$, $r$)}\ncong 0$.
\end{proofs} 

\begin{proposition}\label{arrowsf2}
Let ${\cal C}\hookrightarrow {\cal D}$ be a full embedding of categories such that $\cal D$ is a geometric category with a (strict) initial object $0\in \cal C$ and $\cal C$ is closed in $\cal D$ under pullbacks, cover-mono factorizations and arbitrary unions of subobjects. Given an object $c\in \tilde{{\cal C}}$, a sieve $R$ on $c$ in $\tilde{{\cal C}}$ and a set-indexed collection $\{f_{i}:d_{i}\rightarrow c \textrm{ | } i\in I\}$ of arrows in $\tilde{{\cal C}}$ which are monic in $\cal D$, if $f$ is the union of the subobjects $f_{i}$ (for $i\in I$) in $\Sub_{\cal D}(c)$ then 
\[
\begin{array}{ccl}
f^{\ast}(R)=\emptyset & \textrm{iff} &  \textrm{for each $i\in I$ } f_{i}^{\ast}(R)=\emptyset;\\
f^{\ast}(R) \textrm{ is stably non-empty}  & \textrm{iff} & \textrm{for each $i\in I$ } f_{i}^{\ast}(R) \textrm{ is stably non-empty}.\\  
\end{array}
\]
\end{proposition}

\begin{proofs}
The first part of the proposition follows as an immediate consequence of Proposition \ref{firststep} by using the fact that unions of subobjects in $\cal D$ are stable under pullback. It remains to prove the second part. One implication is obvious, since each $f_{i}$ factors through $f$. To prove the other implication, suppose that for each $i\in I$ $f_{i}^{\ast}(R)$ is stably non-empty. By Propositions \ref{secondstep} and \ref{thirdstep} we can suppose without loss of generality that $R$ is the sieve generated in $\tilde{\cal C}$ by an arrow $r$ in $\tilde{\cal C}$ which is monic in $\cal D$. So we assume that for each $i\in I$ $f_{i}$ is stably joint with $r$ and want to prove that $f$ is stably joint with $r$. Given an arrow $g:e\rightarrow d$ in $\tilde{\cal C}$, let us define for each $i\in I$ $g_{i}$ to be the pullback in $\cal D$ of $g$ along the inclusion $j_{i}:d_{i}\rightarrow c$, as in the following diagram:
\[  
\xymatrix {
y_{i} \ar[d]_{h_{i}} \ar[r]^{g_{i}} & d_{i} \ar[d]^{j_{i}} \ar[rd]^{f_{i}}  & \\
e \ar[r]_{g} & d \ar[r]_{f} & c }
\]\\  
Now, $1_{e}=(f\circ g)^{\ast}(f)=(f\circ g)^{\ast}(\mathbin{\mathop{\textrm{\huge $\cup$}}\limits_{i\in I}}f_{i})=\mathbin{\mathop{\textrm{\huge $\cup$}}\limits_{i\in I}}(f\circ g)^{\ast}(f_{i})=\mathbin{\mathop{\textrm{\huge $\cup$}}\limits_{i\in I}}h_{i}$, hence since $e\ncong 0$ there exists $i\in I$ such that $e_{i}\ncong 0$. So $g_{i}$ is an arrow in $\tilde{\cal C}$ and hence by our assumption $\textrm{p.b.($f_{i}\circ g_{i}$, $r$)}\ncong 0$; the fact that $0$ is a strict initial object then implies that $\textrm{p.b.($f\circ g$, $r$)}\ncong 0$, as required.       
\end{proofs} 
Let us now work under the hypotheses of Corollary \ref{cor2} with the additional assumption that the topology $J$ on $\cal C$ is induced on $\cal C$ by the canonical topology on the geometric category $\cal D$ (equivalently, the $J$-covering sieves are exactly those which contain small covering families). By using such explicit description of the topology $J$, we have that, for each arrow $r:d\rightarrow c$ in $\tilde{\cal C}$ which is monic in $\cal D$, $M_{(r)}\in \tilde{J}(c)$ if and only if there exists a small covering family $\{f_{i} \textrm{ | } i\in I\}$ in $\tilde{\cal C}$ such that for each $i\in I$ either\\
(1) $f_{i}$ is disjoint from $r$\\
or\\
(2) $f_{i}$ is stably joint with $r$.\\
Note that, since $\cal C$ is closed in $\cal D$ under pullbacks and arbitrary unions of subobjects, for each $c\in {\cal C}$, the collection of subobjects in $\Sub_{\cal D}(c)$ which lie (up to isomorphism) in $\cal C$ form a subframe of $\Sub_{\cal D}(c)$; this frame, regarded as a (complete) Heyting algebra, will be denoted by $\Sub^{\cal C}_{\cal D}(c)$.\\  
Suppose that for each $i\in I$ condition (1) holds. Then by Proposition \ref{arrowsf2} there exists a cover $f:dom(f)\epi c$ such that $f$ is disjoint from $r$; this clearly implies ($0$ being strict initial) that $dom(r)=0$, that is $r$ is the zero subobject in the Heyting algebra $\Sub^{\cal C}_{\cal D}(c)$.\\
Suppose that for each $i\in I$ condition (2) holds; again, by Proposition \ref{arrowsf2} we deduce that for every arrow $g:dom(g)\rightarrow c$, $g$ is stably joint with $r$; this is in turn clearly equivalent (by Proposition \ref{arrowsf1}) to saying that $\neg r\cong 0$ in $\Sub^{\cal C}_{\cal D}(c)$.\\
So, provided that $r\ncong 0$ and $\neg r\ncong 0$ in $\Sub^{\cal C}_{\cal D}(c)$, the sets $I_{1}=\{i\in I\textrm{ | $f_{i}$ is disjoint from $r$} \}$ and $I_{2}=\{i\in I\textrm{ | $f_{i}$ is stably joint with $r$} \}$ are both non-empty and we can deduce by Propositions \ref{arrowsf1} and \ref{arrowsf2} that there exists two arrows $f_{1}:dom(f_{1})\rightarrow c$ and $f_{2}:dom(f_{2})\rightarrow c$ in $\tilde{\cal C}$ which are monic in $\cal D$, jointly covering and satisfy: $f_{1}$ is disjoint from $r$ and $f_{2}$ is stably joint from $r$. In terms of the Heyting algebra $\Sub^{\cal C}_{\cal D}(c)$ this condition precisely means that the union of the subobjects $f_{1}$ and $f_{2}$ in $\Sub^{\cal C}_{\cal D}(c)$ is $1_{c}$, $f_{1}\cap r=0$ and $\chi\cap r\neq 0$ for each $\chi$ in $\Sub^{\cal C}_{\cal D}(c)$ such that $\chi\neq0$ and $\chi\leq f_{2}$ (we may suppose - without loss of generality - $\chi$ to be in $\Sub^{\cal C}_{\cal D}(c)$ by Proposition \ref{arrowsf1}). On the other hand, note that the existence of two such arrows $f_{1}$ and $f_{2}$ implies $M_{(r)}\in \tilde{J}(c)$.\\ 
This leads us to introduce the following definition.
\begin{definition}
Let $H$ be an Heyting algebra. Then $H$ is said to satisfy De Morgan property if for each element $r\in H$ such that $r\neq 0$ and $\neg r\neq 0$ there exist elements $f_{1}$ and $f_{2}$ in $H$ satisfying the following conditions:\\
$f_{1}, f_{2}\neq 0$,\\
$f_{1}\vee f_{2}=1$,\\
$f_{1}\wedge r=0$,\\
$\chi \wedge r \neq 0$ for each $\chi\neq 0$ such that $\chi\leq f_{2}$. 
\end{definition}  
\begin{remark}  
This definition can be clearly put also in the following form:\\
an Heyting algebra $H$ satisfies De Morgan property if and only if for each element $r\in H$ such that $r\neq 0$ and $\neg r\neq 0$ there exists a complemented element $f$ in $H$ such that $f\wedge r=0$ and $\chi \wedge r \neq 0$ for each $\chi\neq 0$ such that $\chi\wedge f=0$.\\
\end{remark}
\begin{proposition}
Let $H$ be an Heyting algebra. Then $H$ satisfies De Morgan property if and only if it is a De Morgan algebra.
\end{proposition}     
\begin{proofs}
Let us use the second form of the definition of De Morgan property. In one direction, let us suppose that $H$ satisfies De Morgan property. To prove that $H$ is a De Morgan algebra we need to verify that for each element $r\in H$ we have $\neg r \vee \neg\neg r=1$. Now, if either $r=0$ or $\neg r=0$ this is obvious; if $r\neq0$ and $\neg r\neq0$ then there exists a complemented element $f$ in $H$ such that $f\wedge r=0$ and $\chi \wedge r \neq 0$ for each $\chi\neq 0$ such that $\chi\wedge f=0$. We have that $f\leq \neg r$ as $f\wedge r=0$. To prove that $\neg f\leq \neg\neg r$, observe that $\neg f\wedge \neg r=0$, as otherwise by taking $\chi=\neg f\wedge \neg r$ we would get $\neg f\wedge \neg r\wedge r \neq 0$, a contradiction. So we have $1=f\vee \neg f\leq \neg r \vee \neg\neg r$, that is $\neg r \vee \neg\neg r=1$.\\
Conversely, let us suppose that $H$ is a De Morgan algebra; given $r\in H$ such that $r\neq 0$ and $\neg r\neq 0$, we take $f$ to be the complemented element $\neg r$; this does the work because obviously $f\wedge r=0$ and given $\chi\neq 0$ such that $\chi\wedge f=0$, $\chi \wedge r \neq 0$ as otherwise we would have $\chi\leq \neg r$ and hence $\chi\leq \neg r \wedge \neg\neg r=0$.   
\end{proofs}    
So we have arrived at the following result.
\begin{theorem}\label{thm}
Let $({\cal C}, J)$ be a Grothendieck site and ${\cal C}\hookrightarrow {\cal D}$ be a full embedding of categories such that $\cal D$ is a geometric category with a (strict) initial object $0\in \cal C$ and $\cal C$ is closed in $\cal D$ under pullbacks, cover-mono factorizations and arbitrary unions of subobjects. If $J$ is the Grothendieck topology on $\cal C$ induced by the canonical topology on $\cal D$ then $\Sh({\cal C}, J)$ is a De Morgan topos if and only if for each object $c\in \cal C$ the Heyting algebra $\Sub^{\cal C}_{\cal D}(c)$ satisfies De Morgan property (equivalently, is a De Morgan algebra).
\end{theorem}\qed
Finally, let us consider how our simplification method can be adapted to the Boolean case.\\
From the considerations in the first section we deduce the following criterion: provided that every $J$-covering sieve is non-empty, $\Sh({\cal C}, J)$ is Boolean if and only if every stably non-empty sieve in $\tilde{\cal C}$ is $\tilde{J}$-covering.
Now, if $\cal D$ is a regular (resp. geometric) category and $J$ is the Grothendieck topology on $\cal C$ induced via the embedding ${\cal C}\hookrightarrow {\cal D}$ by the canonical topology on $\cal D$, Proposition \ref{secondstep} (resp. Proposition \ref{thirdstep}) enables us to restrict our attention to sieves $R$ which are generated by a family of arrows which are monic in $\cal D$ (resp. by a single arrow which is monic in $\cal D$), as in Corollary \ref{cor1} (resp. Corollary \ref{cor2}). In fact, the following results hold.
\begin{corollary}\label{cor1b}
Let $({\cal C}, J)$ be a Grothendieck site and ${\cal C}\hookrightarrow {\cal D}$ a full embedding of categories such that $\cal D$ is a regular category having a strict initial object $0\in \cal C$ and $\cal C$ is closed in $\cal D$ under pullbacks and cover-mono factorizations. Then\\
(a) If $J$ is the Grothendieck topology on $\cal C$ induced by the canonical topology on $\cal D$ then $\Sh({\cal C}, J)$ is Boolean if and only if every stably non-empty sieve $R$ in $\tilde{\cal C}$ generated in $\tilde{\cal C}$ by morphisms which are monic in $\cal D$ is a $\tilde{J}$-covering sieve (equivalently, the maximal sieve).\\
(b) If $\cal D$ is a coherent category and $J$ is the Grothendieck topology on $\cal C$ induced by the canonical topology on $\cal D$ then $\Sh({\cal C}, J)$ is Boolean if and only if every stably non-empty sieve $R$ in $\tilde{\cal C}$ generated in $\tilde{\cal C}$ by morphisms which are monic in $\cal D$ contains a finite covering family.\\    
\end{corollary}
\begin{proofs}
(a) One direction is obvious. Let us prove the other one. If $R$ is a stably non-empty sieve in $\tilde{\cal C}$ on an object $c$ then by Proposition \ref{secondstep} the sieve $R'$ generated by the images in $\cal D$ of the morphisms in $R$ is stably non-empty and hence $\tilde{J}$-covering. Then, $\tilde{J}$ being induced by the canonical topology on the regular category $\cal D$, there exists a morphism in $\tilde{\cal C}$ which is a cover in $\cal D$ and belongs to $R'$. Thus, the identity $1_{c}$ factors through one of the generating morphisms of $R'$, that is there exists a morphism in $R$ whose image is isomorphic to the identity, i.e. which is a cover in $\cal D$; hence $R$ is a $\tilde{J}$-covering sieve.\\
(b) One direction is obvious. In the other direction, given a stably non-empty sieve $R$, consider the sieve $R'$ as above. Then $R'$ is $\tilde{J}$-covering and hence, $\tilde{J}$ being induced by the canonical topology on the coherent category $\cal D$, $R'$ contains a finite covering family. In fact, $R'$ being generated by monomorphisms, we may clearly suppose the members of such a family to belong to this collection of monomorphisms. Then $R$ contains a finite covering family (take the arrows whose images are in the family above), and hence is $\tilde{J}$-covering.    
     
\end{proofs}

\begin{corollary}\label{cor2b}
Let $({\cal C}, J)$ be a Grothendieck site and ${\cal C}\hookrightarrow {\cal D}$ a full embedding of categories such that $\cal D$ is a geometric category with a (strict) initial object $0\in \cal C$ and $\cal C$ is closed in $\cal D$ under pullbacks, cover-mono factorizations and arbitrary unions of subobjects of objects in $\cal C$. If $J$ is the Grothendieck topology on $\cal C$ induced by the canonical topology on $\cal D$ then $\Sh({\cal C}, J)$ is Boolean if and only if for each arrow $r$ in $\tilde{\cal C}$ such that $r$ is monic in $\cal D$ and $(r)$ is stably non-empty, $(r)$ is a $\tilde{J}$-covering sieve (equivalently, the maximal sieve). 
\end{corollary}
\begin{proofs}
One direction being obvious, let us prove the other one. By Corollary \ref{cor1b} we can restrict our attention to sieves $R$ generated by a set $\{r_{i} \textrm{ | } i\in I \}$ of arrows in $\tilde{\cal C}$ which are monomorphisms in $\cal D$. Denoted by $r$ the union of these monomorphisms in $\Sub_{\cal D}(c)$, we have by Proposition \ref{thirdstep} that $(r)$ is stably non-empty. Then $(r)$ is a $\tilde{J}$-covering sieve, that is ($\tilde{J}$ being induced by the canonical topology on the geometric category $\cal D$) $(r)$ contains a small covering family of arrows lying in $\tilde{\cal C}$. Now, the fact that $r$ is monic in $\cal D$ implies that the sieve $(r)$ is closed in $\cal D$ under taking images and unions of subobjects in $\tilde{\cal C}$, so $1_{c}\in (r)$. Thus, $\{r_{i} \textrm{ | } i\in I \}$ is a small covering family and hence $R$ is a $\tilde{J}$-covering sieve.   
\end{proofs}

Analogously to the De Morgan case, we are led to introduce the following notion.
\begin{definition}
Let $H$ be an Heyting algebra. Then $H$ is said to satisfy the Boolean property if the only element $r\in H$ such that for each $\chi\in H$, $\chi\neq 0$ implies $\chi\wedge r\neq 0$, is $1$. 
\end{definition}
  
\begin{proposition}
Let $H$ be an Heyting algebra. Then $H$ satisfies the Boolean property if and only if it is a Boolean algebra.
\end{proposition}     
\begin{proofs}
In one direction, let us suppose that $H$ satisfies the Boolean property. To prove that $H$ is a Boolean algebra we need to verify that for each element $r\in H$ we have $r \vee \neg r=1$. Now, for any $\chi\neq 0$ we have $\chi\wedge (r\vee \neg r)=(\chi\wedge r)\vee (\chi\wedge \neg r)\neq 0$ because otherwise we would have $\chi\wedge r=0$, $\chi\wedge \neg r=0$ and hence $\chi \leq \neg r \wedge \neg\neg r=0$, which is absurd; so we conclude that $r\vee \neg r=1$.\\
Conversely, let us suppose that $H$ is a Boolean algebra; then given $r\in H$ such that $r\neq 1$, $\chi=\neg r$ satisfies $\chi\neq 0$ and $\chi \wedge r=0$.
\end{proofs}  
The analogue of Theorem \ref{thm} is then given by the following result.
\begin{theorem}\label{thmb}
Let $({\cal C}, J)$ be a Grothendieck site and ${\cal C}\hookrightarrow {\cal D}$ be a full embedding of categories such that $\cal D$ is a geometric category with a (strict) initial object $0\in \cal C$ and $\cal C$ is closed in $\cal D$ under pullbacks, cover-mono factorizations and arbitrary unions of subobjects. If $J$ is the Grothendieck topology on $\cal C$ induced by the canonical topology on $\cal D$ then $\Sh({\cal C}, J)$ is Boolean if and only if for each object $c\in \cal C$ the Heyting algebra $\Sub^{\cal C}_{\cal D}(c)$ satisfies Boolean property (equivalently, is a Boolean algebra).
\end{theorem}
\begin{proofs}
This immediately follows from Corollary \ref{cor2b} by using Proposition \ref{arrowsf1}.  
\end{proofs}

\section{Applications}
\subsection{Syntactic criteria}
Given a geometric theory $\mathbb T$, we say that $\mathbb T$ is a De Morgan (resp. Boolean) theory if its classifying topos $\Set[\mathbb T]$ satisfies De Morgan's law (resp. is Boolean).\\
In this section we show how it is possible to deduce from Theorem \ref{thm} (resp. Theorem \ref{thmb}) in the last section a syntactic criterion for a geometric theory to be a De Morgan (resp. Boolean) theory.\\
We recall from \cite{El2} that, given a geometric theory $\mathbb T$, its classifying topos $\Set[{\mathbb T}]$ for $\mathbb T$ can be represented as the category $\Sh({\cal C}_{\mathbb T}, J_{\mathbb T})$ of sheaves on the syntactic category ${\cal C}_{\mathbb T}$ of $\mathbb T$ with respect to the syntactic topology $J_{\mathbb T}$ on it (i.e. the canonical topology on the geometric category ${\cal C}_{\mathbb T}$). Taking ${\cal C}={\cal D}={\cal C}_{\mathbb T}$ the hypotheses of Theorem \ref{thm} (resp. Theorem \ref{thmb}) are clearly satisfied, so we obtain the following criterion: given a geometric theory $\mathbb T$, $\mathbb T$ is a De Morgan (resp. Boolean) theory if and only if the subobject lattices $\Sub_{{\cal C}_{\mathbb T}}(c)$ (for $c\in {\cal C}_{\mathbb T}$) are all De Morgan (resp. Boolean) algebras (equivalently, they satisfy De Morgan (resp. Boolean) property). In fact, it is possible to rephrase this latter condition as a syntactic property of the geometric theory $\mathbb T$, as in the following results.\\
Below, in the context of a geometric theory $\mathbb T$ over a signature $\Sigma$, a geometric formula $\phi(\vec{x})$ is said to be consistent if the sequent $\phi(\vec{x}) \: \vdash_{\vec{x}}\: \bot$ is not provable in $\mathbb T$.             
\begin{theorem}\label{teosyint}
Let $\mathbb T$ be a geometric theory over a signature $\Sigma$. Then $\mathbb T$ is a De Morgan theory if and only if for every consistent geometric formula $\phi(\vec{x})$ over $\Sigma$ such that $\top \:\vdash_{\vec{x}}\: \phi(\vec{x})$ is not provable in $\mathbb T$, there exists two consistent geometric formulae $\psi_{1}(\vec{x})$ and $\psi_{2}(\vec{x})$ over $\Sigma$ in the same context such that:\\
$\top \:\vdash_{\vec{x}}\: \psi_{1}(\vec{x})\vee \psi_{2}(\vec{x})$ is provable in $\mathbb T$,\\
$\psi_{1}(\vec{x})\wedge \phi(\vec{x}) \:\vdash_{\vec{x}}\: \bot$ is provable in $\mathbb T$ and\\
$\chi(\vec{x}) \wedge \phi(\vec{x})$ is consistent for every consistent geometric formula $\chi(\vec{x})$ over $\Sigma$ in the same context such that $\chi(\vec{x}) \:\vdash_{\vec{x}}\: \psi_{2}(\vec{x})$ is provable in $\mathbb T$.\\ 
\end{theorem}
\begin{proofs}
For each geometric formula $\phi(\vec{x})$, there is an obvious monomorphism $\{\vec{x}.\phi(\vec{x})\}\rightarrow \{\vec{x}.\top\}$ in the syntactic category ${\cal C}_{\mathbb T}$, so by Remark \ref{rmks}(b), we can restrict our attention to the subobject lattices $\Sub_{{\cal C}_{\mathbb T}}(\{\vec{x}.\top\})$. Our thesis then follows from the explicit description of subobjects in the syntactic category ${\cal C}_{\mathbb T}$ given by Lemma D1.4.4(iv) \cite{El2}.  
\end{proofs}
\begin{theorem}\label{teosyintb}
Let $\mathbb T$ be a geometric theory over a signature $\Sigma$. Then $\mathbb T$ is a Boolean theory if and only if every consistent geometric formula $\phi(\vec{x})$ over $\Sigma$ such that $\chi(\vec{x}) \wedge \phi(\vec{x})$ is consistent for each consistent geometric formula $\chi(\vec{x})$ over $\Sigma$ in the same context is provable equivalent to $\top$ in $\mathbb T$; equivalently, for every geometric formula $\phi(\vec{x})$ over $\Sigma$ there is a geometric formula $\chi(\vec{x})$ over $\Sigma$ in the same context such that $\phi(\vec{x}) \wedge \chi(\vec{x}) \:\vdash_{\vec{x}}\: \bot$ and $\top \:\vdash_{\vec{x}}\: \phi(\vec{x}) \vee \chi(\vec{x})$ are provable in $\mathbb T$.   
\end{theorem}
\begin{proofs}
Similar to the proof of Theorem \ref{teosyint}.
\end{proofs}
Now, let us suppose that $\mathbb T$ is a coherent theory over a signature $\Sigma$. The classifying topos $\Set[\mathbb T]$ can be represented as the category $\Sh({\cal C}^{\textrm{coh}}_{\mathbb T}, J^{\textrm{coh}}_{\mathbb T})$ of sheaves on the coherent syntactic category ${\cal C}^{\textrm{coh}}_{\mathbb T}$ of $\mathbb T$ with respect to the canonical topology $J^{\textrm{coh}}_{\mathbb T}$ on it, that is the topology having as covering sieves those which contain finite covering families.\\
From Corollary \ref{cor1} and Proposition \ref{arrowsf1}, by arguing as above, we immediately obtain the following result.
\begin{theorem}\label{teosyintcoh}
Let $\mathbb T$ be a coherent theory over a signature $\Sigma$. Then $\mathbb T$ is a De Morgan theory if and only if for every family $\{\phi_{i}(\vec{x}) \textrm{ | } i\in I\}$ of consistent coherent formulae over $\Sigma$ in the same context there exists a finite family $\psi_{1}(\vec{x}), \psi_{2}(\vec{x}), \ldots, \psi_{n}(\vec{x})$ of consistent coherent formulae over $\Sigma$ in the same context such that:\\
$\top \:\vdash_{\vec{x}}\: \mathbin{\mathop{\textrm{\huge $\vee$}}\limits_{1\leq j\leq n}}\psi_{j}(\vec{x})$ is provable in $\mathbb T$ and for each $1\leq j\leq n$\\
either $\psi_{j}(\vec{x})\wedge \phi_{i}(\vec{x}) \:\vdash_{\vec{x}}\: \bot$ is provable in $\mathbb T$ for all $i\in I$ or\\
for every consistent geometric formula $\chi(\vec{x})$ such that $\chi(\vec{x}) \:\vdash_{\vec{x}}\: \psi_{j}(\vec{x})$ is provable in $\mathbb T$ there exists $i\in I$ such that $\chi(\vec{x}) \wedge \phi_{i}(\vec{x})$ is consistent.\\ 
\end{theorem}\qed

\begin{theorem}\label{teosyintbcoh}
Let $\mathbb T$ be a coherent theory over a signature $\Sigma$. Then $\mathbb T$ is a Boolean theory if and only if for every family $\{\phi_{i}(\vec{x}) \textrm{ | } i\in I\}$ of consistent coherent formulae over $\Sigma$ in the same context with the property that for each coherent consistent formula $\chi(\vec{x})$ over $\Sigma$ in the same context there exists $i\in I$ such that $\phi_{i}(\vec{x})\wedge \chi(\vec{x})$ is consistent, there exists a finite subset $J\subseteq I$ such that\\  
$\top \:\vdash_{\vec{x}}\: \mathbin{\mathop{\textrm{\huge $\vee$}}\limits_{j\in J}}\phi_{j}(\vec{x})$ is provable in $\mathbb T$.\\
\end{theorem}
\begin{proofs}
This follows as an immediate consequence of Corollary \ref{cor1b}(b) by identifying formulas with the corresponding monomorphisms in the coherent syntactic category (as in the proof of Theorem \ref{teosyint}).
\end{proofs}

Note that for families $\{\phi_{i}(\vec{x}) \textrm{ | } i\in I\}$ formed by just one element, Theorem \ref{teosyintbcoh} says precisely that the subobject lattices in the coherent syntactic category ${\cal C}^{\textrm{coh}}_{\mathbb T}$ satisfy Boolean property (equivalently, are Boolean algebras). If $\mathbb T$ is Boolean, this also implies that they are finite (cfr. the proof of Theorem D3.4.3 \cite{El2}). On the other hand, given a coherent theory $\mathbb T$ such that all the subobject lattices in ${\cal C}^{\textrm{coh}}_{\mathbb T}$ are finite, we may immediately deduce from Theorem \ref{teosyintbcoh} that if they are all also Boolean algebras, $\mathbb T$ is Boolean (under these hypotheses, all the families $\{\phi_{i}(\vec{x}) \textrm{ | } i\in I\}$ in the statement of the theorem are finite and hence they can be replaced - for our purposes - by the singletons $\{\phi(\vec{x})\}$, where $\phi(\vec{x})$ is the finite disjunction of all the $\phi_{i}(\vec{x})$). In view of Remark \ref{rmks}(b), this proves that a coherent theory $\mathbb T$ is Boolean if and only if all the subobject lattices in ${\cal C}^{\textrm{coh}}_{\mathbb T}$ of the form $\Sub(\{\vec{x}.\top\})$ are finite Boolean algebras; we note that this is essentially the content of Theorem D3.4.6 \cite{El2}.\\  
Suppose now that $\mathbb T$ is a regular theory. The classifying topos $\Set[\mathbb T]$ can be represented as the category $\Sh({\cal C}^{\textrm{reg}}_{\mathbb T}, J_{\mathbb T}^{\textrm{reg}})$ of sheaves on the regular syntactic category ${\cal C}^{\textrm{reg}}_{\mathbb T}$ of $\mathbb T$ with respect to the canonical topology $J_{\mathbb T}^{\textrm{reg}}$ on it, that is the topology having as covering sieves those which contain a cover. Since the category ${\cal C}^{\textrm{reg}}_{\mathbb T}$ satisfies the right Ore condition (being cartesian) and the topology $J_{\mathbb T}^{\textrm{reg}}$ has no empty covering sieves, we deduce from Corollary \ref{corfond} that $\Set[\mathbb T]\simeq \Sh({\cal C}^{\textrm{reg}}_{\mathbb T}, J_{\mathbb T}^{\textrm{reg}})$ is a De Morgan topos. We have thus proved the following result.

\begin{theorem}\label{teosyintreg}
Let $\mathbb T$ be a regular theory. Then $\mathbb T$ is a De Morgan theory.
\end{theorem}\qed  
Concerning the Boolean case, we have the following result.

\begin{theorem}\label{teosyintbreg}
Let $\mathbb T$ be a regular theory over a signature $\Sigma$. Then $\mathbb T$ is a Boolean theory if and only if for every family $\{\phi_{i}(\vec{x}) \textrm{ | } i\in I\}$ of consistent regular formulae over $\Sigma$ in the same context such that for each regular consistent formula $\chi(\vec{x})$ over $\Sigma$ in the same context there exists $i\in I$ such that $\phi_{i}(\vec{x})\wedge \chi(\vec{x})$ is consistent, there exists $i\in I$ such that $\top \: \vdash_{\vec{x}}\: \phi_{i}(\vec{x})$ is provable in $\mathbb T$.
\end{theorem}
\begin{proofs}
This immediately follows from Corollary \ref{cor1b}(a) (by the usual identification of formulas with monomorphisms in the relevant syntactic category). 
\end{proofs}

\subsection{Separating sets for Grothendieck toposes}
We observe that our simplification method can be easily modified to obtain a version of it for $\infty$-pretoposes in place of geometric categories; in particular, we have the following result.
\begin{theorem}\label{thmpretopos}
Let $({\cal C}, J)$ be a Grothendieck site and ${\cal C}\hookrightarrow {\cal D}$ be a full embedding of categories such that $\cal D$ is an $\infty$-pretopos $\cal D$ with a (strict) initial object $0\in \cal C$ and $\cal C$ is closed in $\cal D$ under pullbacks, cover-mono factorizations and arbitrary unions of subobjects. If $J$ is the Grothendieck topology on $\cal C$ induced by the canonical topology on $\cal D$ then $\Sh({\cal C}, J)$ is a De Morgan topos (resp. a Boolean topos) if and only if for each object $c\in \cal C$ the subobject lattice $\Sub^{\cal C}_{\cal D}(c)$ is a De Morgan algebra (resp. a Boolean algebra). 
\end{theorem}\qed
From this theorem one can immediately deduce that if ${\cal C}\hookrightarrow {\cal E}$ is a separating set for a Grothendieck topos $\cal E$ which is closed under in $\cal E$ under pullbacks and under taking subobjects in $\cal E$ then $\cal E$ is a De Morgan (resp. Boolean) topos if and only if for each $c\in \cal C$, $\Sub_{\cal E}(c)$ is a De Morgan (resp. Boolean) algebra. In fact, the following more general result hold.
\begin{theorem}\label{separsets}
Let $\cal E$ be an $\infty$-pretopos with a separating set $\cal C$. Then $\cal E$ is a De Morgan (resp. Boolean) topos if and only if for each $c\in \cal C$, $\Sub_{\cal E}(c)$ is a De Morgan (resp. Boolean) algebra.  
\end{theorem}
\begin{proofs}
Given an elementary topos $\cal E$, it is well-known that $\cal E$ is a De Morgan (resp. Boolean) topos if and only if all the subobject lattices $\Sub_{\cal E}(e)$ for $e\in \cal E$ are De Morgan (resp. Boolean) algebras; here we want to show that, under our hypotheses, it is enough to check that all the subobject lattices $\Sub_{\cal E}(c)$ for $c\in \cal C$ are. Given a Grothendieck topos $\cal E$ and an object $e\in \cal E$, if $\cal C$ is a separating set for $\cal E$ then $e$ can be expressed as a quotient of a coproduct of objects in $\cal C$, that is there exists a family $\{c_{i} \textrm{ | } i\in I\}$ (where $I$ is a set) of objects in $\cal C$ and an epimorphism $p:\coprod_{i\in I}c_{i}\epi e$. Since $p$ is an epimorphism, then the pullback functor $p^{\ast}:\Sub_{\cal E}(e) \rightarrow \Sub_{\cal E}(\coprod_{i\in I}c_{i})\cong \prod_{i\in I}\Sub_{\cal E}(c_{i})$ is (logical and) conservative (cfr. Example 4.2.7(a) p. 181 \cite{El}), hence $\Sub_{\cal E}(e)$ is a De Morgan (resp. Boolean) algebra if all the $\Sub_{\cal E}(c_{i})$ are. 
\end{proofs}

\begin{remarks}
(a) The classifying topos $\Set[\mathbb T]$ of a geometric theory $\mathbb T$ is the $\infty$-pretopos generated by the geometric syntactic category ${\cal C}_{\mathbb T}$ of $\mathbb T$ (cfr. Proposition D3.1.12 \cite{El2}), so the objects of ${\cal C}_{\mathbb T}$ form a separating set for $\Set[\mathbb T]$; hence the hypotheses of Theorem \ref{separsets} are satisfied and we get Theorems \ref{teosyint} and \ref{teosyintb} as an application.\\
(b) Another case of interest in which the theorem can be applied is when we have a pre-bound $B$ for $\cal E$; indeed, the subobjects of finite powers $B^{n}$ form a separating set for $\cal E$ and hence, in view of Remark \ref{rmks}(b), we obtain the following characterization: $\cal E$ is a De Morgan (resp. Boolean) topos if and only if for each natural number $n$ the Heyting algebra $\Sub(B^{n})$ is a De Morgan (resp. Boolean) algebra; in particular, if $\Set[\mathbb T]$ is the classifying topos of a one-sorted geometric theory $\mathbb T$ then the underlying object $M_{\mathbb T}\in \Set[\mathbb T]$ of the universal $\mathbb T$-model is a pre-bound for $\Set[\mathbb T]$ so we obtain the following criterion: $\mathbb T$ is a De Morgan (resp. Boolean) theory if and only if all the lattices $\Sub({M_{\mathbb T}}^{n})$ (for $n$ natural number) are De Morgan (resp. Boolean) algebras.   
\end{remarks} 

\subsection{Topological interpretations}
In section $2$, we have introduced the notions of an Heying algebra satisfying De Morgan (resp. Boolean) property; in this section, we show that these notions have a clear topological meaning in terms of the (generalized) locale corresponding to the Heyting algebra.\\
Given an Heyting algebra $H$, we can consider it as a generalized locale (recall that locales are the same thing as \emph{complete} Heyting algebras). We recall that any sublocale $Y$ of a given locale $X$ has a closure $\overline{Y}$ ; in particular, if $Y$ is the open sublocale of a locale $H$ corresponding to an element $a\in H$, then $\overline{Y}$ is the closed sublocale of $H$ corresponding to the element $\neg a\in H$ (note that this notion of closure of an open sublocale is liftable from the context of locales to that of generalized locales). Now, by using the characterization of the closure of a sublocale given in the proof of Proposition \ref{loctop}, it is easy to see that De Morgan property on an Heyting algebra is equivalent to the statement that the corresponding (generalized) locale is extremally disconnected (i.e. the closure of any open sublocale is open), while the Boolean property is equivalent to saying that the (generalized) locale is almost discrete, (i.e. the only non-empty dense open sublocale is the whole locale, in other words every open sublocale is closed).\\
By regarding a locale $L$ as a geometric category, Theorems \ref{thm} and \ref{thmb} (together with Remark \ref{rmks}(b)) then give the following results.
\begin{theorem}
Let $L$ be a locale. Then $\Sh(L)$ is a De Morgan topos if and only if $L$ is a De Morgan algebra (equivalently, satisfies De Morgan property), if and only if $L$ is extremally disconnected. 
\end{theorem}\qed
\begin{theorem}
Let $L$ be a locale. Then $\Sh(L)$ is a Boolean topos if and only if $L$ is a Boolean algebra (equivalently, satisfies Boolean property), if and only if $L$ is almost discrete. 
\end{theorem}\qed
These results are more or less well-known, but we feel that, by introducing the concepts of De Morgan and Boolean property, we have added another viewpoint that clarifies the interplay between the topological and logical notions.\\
Now, let $\mathbb T$ be a geometric theory over a signature $\Sigma$ and $M$ a $\mathbb T$-model in a Grothendieck topos $\cal E$. For each context $\vec{x}=(x_{1}^{A_{1}}, \ldots, x_{n}^{A_{n}})$ over $\Sigma$, the subobjects of $MA_{1}\times \ldots \times MA_{n}$ of the form $[[\phi(\vec{x})]]$, where $\phi(\vec{x})$ is a geometric formula in the context $\vec{x}$ over $\Sigma$ (here $[[\vec{x}. \phi]]$ denotes the interpretation of the formula $\phi(\vec{x})$ in the model $M$) clearly form a subframe of $\Sub_{\cal E}(MA_{1}\times \ldots \times MA_{n})$. The locale corresponding to this frame will be denoted by $\textrm{Def}_{\vec{x}}^{\textrm{geom}}(M)$; we note that, at least when $\cal E$ is the topos $\Set$, this locale is spatial, so that we have a topological space $\textrm{Def}_{\vec{x}}^{\textrm{geom}}(M)$ whose open subsets are exactly the subsets of $MA_{1}\times \ldots \times MA_{n}$ which are definable by geometric formulas.\\
Given a geometric theory $\mathbb T$, we now consider how the property of the classifying topos $\Set[\mathbb T]$ to be De Morgan (or Boolean) reflects into topological or logical properties of the locales $\textrm{Def}_{\vec{x}}^{\textrm{geom}}(M)$. First, we note that there is a geometric surjective functor $\textrm{Int}_{\vec{x}}^{M}:\Sub_{{\cal C}^{\textrm{geom}}_{\mathbb T}}(\{\vec{x}.\top\})\rightarrow \textrm{Def}_{\vec{x}}^{\textrm{geom}}(M)$ which sends each formula $\phi(\vec{x})$ (identified with the corresponding subobject $\{\vec{x}.\phi(\vec{x})\}\rightarrow \{\vec{x}.\top\}$ in ${\cal C}^{\textrm{geom}}_{\mathbb T}$) to the interpretation $[[\vec{x}. \phi]]$ in the model $M$. As a consequence of the fact that $\textrm{Int}_{\vec{x}}^{M}$ is geometric we deduce that if $\Sub_{{\cal C}^{\textrm{geom}}_{\mathbb T}}(\{\vec{x}.\top\})$ is a Boolean algebra then $\textrm{Def}_{\vec{x}}^{\textrm{geom}}(M)$ is also a Boolean algebra; however, it is not true in general that if $\Sub_{{\cal C}^{\textrm{geom}}_{\mathbb T}}(\{\vec{x}.\top\})$ is a De Morgan algebra then $\textrm{Def}_{\vec{x}}^{\textrm{geom}}(M)$ is a De Morgan algebra. If $M$ is a conservative $\mathbb T$-model, then clearly $\textrm{Int}_{\vec{x}}^{M}$ is conservative and hence an isomorphism, so $\Sub_{{\cal C}^{\textrm{geom}}_{\mathbb T}}(\{\vec{x}.\top\})$ is a De Morgan (resp. Boolean) algebra if and only if $\textrm{Def}_{\vec{x}}^{\textrm{geom}}(M)$ is. As an application of this, consider the universal model $M_{\mathbb T}$ of a geometric theory $\mathbb T$ lying in the classifying topos $\Set[\mathbb T]$; in view of our characterization saying that $\Set[\mathbb T]$ is De Morgan (resp. Boolean) if and only if all the subobject lattices of the form so $\Sub_{{\cal C}^{\textrm{geom}}_{\mathbb T}}(\{\vec{x}.\top\})$ are De Morgan (resp. Boolean) algebras, we obtain the following criterion: $\mathbb T$ is a De Morgan (resp. Boolean) theory if and only if all the $\textrm{Def}_{\vec{x}}^{\textrm{geom}}(M_{\mathbb T})$ are De Morgan (resp. Boolean) algebras.

\section{Model-theoretic characterizations}
Let us first introduce some notation. Given two Grothendieck topologies $J$ and $J'$ on a given category $\cal C$, we write $J'\subseteq J$ to mean that every $J'$-covering sieve is a $J$-covering sieve (equivalently, $J'\leq J$ as topologies on the topos $[{\cal C}^{\textrm{op}}, \Set]$). Given a Grothendieck topology $J$ on a category $\cal C$, we denote by $a_{J}:[{\cal C}^{\textrm{op}}, \Set]\rightarrow \Sh({\cal C}, J)$ the associated sheaf functor.\\
Given a Grothendieck topos $\cal E$ and a category $\cal C$, we write $\bf Flat(\cal{C}, \cal{E})$ for the category of flat functors $\st{\cal{C}}{\cal{E}}$ and natural transformations between them; for a Grothendieck topology $J$ on $\cal C$, ${\bf Flat}_{J}(\cal{C}, \cal{E})$ will denote the full subcategory of $J$-continuous flat functors $\st{\cal{C}}{\cal{E}}$. The $2$-category of Grothendieck toposes, geometric morphisms and geometric transformations between them will be denoted by $\mathfrak{BTop}$ and, given two Grothendieck toposes $\cal E$ and $\cal F$, we will write ${\bf Geom}({\cal E},{\cal F})$ for the category of geometric morphisms $\st{{\cal E}}{{\cal F}}$ and geometric transformations between them.\\    

\begin{lemma}\label{lemma1}
Let $J$ and $J'$ be two Grothendieck topologies on a given category $\cal C$. Then $J'\subseteq J$ if and only if for each Grothendieck topos $\cal E$ every $J$-continuous flat functor $\st{\cal{C}}{\cal{E}}$ is $J'$-continuous, equivalently the functor $a_{J}\circ y:{\cal C}\rightarrow \Sh({\cal C}, J)$ is $J'$-continuous (where $y:{\cal C}\rightarrow [{\cal C}^{\textrm{op}}, \Set]$ is the Yoneda embedding).
\end{lemma}
\begin{proofs}
We recall that there is an equivalence of categories ${\bf Flat}_{J}({\cal C}, {\cal E})\simeq {\bf Geom}({\cal E},\Sh({\cal C}, J))$, which is natural in ${\cal E}\in \mathfrak{BTop}$. By this equivalence, requiring that for each Grothendieck topos $\cal E$ there is an inclusion ${\bf Flat}_{J}({\cal C}, {\cal E})\subseteq {\bf Flat}_{J'}({\cal C}, {\cal E})$ as in the statement of the lemma, is equivalent to demanding that for each ${\cal E}\in \mathfrak{BTop}$ there is a commutative diagram  
\[  
\xymatrix {
{\bf Geom}({\cal E},\Sh({\cal C}, J)) \ar[dr] \ar[rr]  &  & {\bf Geom}({\cal E},\Sh({\cal C}, J')) \ar[dl]  \\
& {\bf Geom}({\cal E},[{\cal C}^{\textrm{op}}, \Set]) &}
\]\\  
which is natural in ${\cal E}\in \mathfrak{BTop}$, where the two diagonal arrows are the obvious ones induced by the inclusions.
This is in turn equivalent, by Yoneda, to requiring that the geometric inclusion $\Sh({\cal C}, J)\hookrightarrow [{\cal C}^{\textrm{op}}, \Set]$ factors through the inclusion $\Sh({\cal C}, J')\hookrightarrow [{\cal C}^{\textrm{op}}, \Set]$ (equivalently the flat $J$-continous functor $a_{J}\circ y:{\cal C}\rightarrow \Sh({\cal C}, J)$ is $J'$-continuous); and from the theory of elementary toposes we know that this happens precisely when $J'\subseteq J$. 
\end{proofs}

Given a Grothendieck site $({\cal C}, J)$, let us consider the ``reduced'' site $(\tilde{{\cal C}},\tilde{J})$, as in Section $2$. As we have already remarked, there is an equivalence of categories $\Sh({\cal C}, J)\simeq \Sh(\tilde{{\cal C}},\tilde{J})$, given by the Comparison Lemma. This equivalence is in fact a geometric equivalence of toposes $\tau: \Sh(\tilde{{\cal C}},\tilde{J})\rightarrow \Sh({\cal C}, J)$, having as its inverse image the obvious restriction functor. Indeed, from the proof of Theorem C2.2.3 \cite{El2} we see (by invoking the uniqueness - up to isomorphism - of right adjoints), that the geometric morphism $l:[\tilde{{\cal C}}^{\textrm{op}}, \Set]\rightarrow [{\cal C}^{\textrm{op}}, \Set]$ induced by the inclusion $\tilde{\cal C}^{\textrm{op}}\hookrightarrow {\cal C}^{\textrm{op}}$ (as in Example A4.1.4 \cite{El}), restricts to the equivalence $\tau$ between the subtoposes $\Sh(\tilde{{\cal C}},\tilde{J})$ and $\Sh({\cal C}, J)$, that is we have a commutative diagram 
\[  
\xymatrix {
\Sh(\tilde{{\cal C}}, \tilde{J}) \ar[d]_{\tau} \ar[r]  & [\tilde{{\cal C}}^{\textrm{op}}, \Set]  \ar[d]^{l}  \\
\Sh({\cal C}, J) \ar[r] & [{\cal C}^{\textrm{op}}, \Set]  }
\]
in $\mathfrak{BTop}$ (where the horizontal arrows are the obvious geometric inclusions). 
Now, for each Grothendieck topos $\cal E$, the equivalence of categories $-\circ \tau: {\bf Geom}({\cal E},\Sh({\cal C}, J)) \rightarrow {\bf Geom}({\cal E},\Sh(\tilde{{\cal C}},\tilde{J}))$ obtained by composing with $\tau$, induce, via the equivalences ${\bf Geom}({\cal E},\Sh({\cal C}, J)) \simeq {\bf Flat}_{J}({\cal C}, {\cal E})$ and ${\bf Geom}({\cal E},\Sh(\tilde{{\cal C}}, \tilde{J})) \simeq {\bf Flat}_{\tilde{J}}(\tilde{{\cal C}}, {\cal E})$, an equivalence of categories ${\bf Flat}_{J}({\cal C}, {\cal E})\simeq {\bf Flat}_{\tilde{J}}(\tilde{{\cal C}}, {\cal E})$, whose explicit description is given by the following lemma.

\begin{lemma}\label{lemma2}
With the above notation, the equivalence ${\bf Flat}_{J}({\cal C}, {\cal E})\simeq {\bf Flat}_{\tilde{J}}(\tilde{{\cal C}},  {\cal E})$ has the following description: one half of the equivalence sends a $J$-continuous flat functor $F:{\cal C}\rightarrow {\cal E}$ to its restriction $F\textrm{|}_{\tilde{\cal C}}:\tilde{\cal C}\rightarrow {\cal E}$ to the category $\tilde{\cal C}$, while the other half of the equivalence sends a $\tilde{J}$-continous flat functor $G:\tilde{\cal C}\rightarrow {\cal E}$ to the its extension $\overline{G}:{\cal C}\rightarrow {\cal E}$ to $\cal C$ obtained by putting $\overline{G}(c)=0$ for each $c\in \cal C$ not in $\tilde{\cal C}$.
\end{lemma}
\begin{proofs}
The equivalence ${\bf Geom}({\cal E},\Sh({\cal C}, J)) \simeq {\bf Flat}_{J}({\cal C}, {\cal E})$ (resp. ${\bf Geom}({\cal E},\Sh(\tilde{{\cal C}}, \tilde{J})) \simeq {\bf Flat}_{\tilde{J}}(\tilde{{\cal C}}, {\cal E})$), sends a geometric morphism $f:{\cal E}\rightarrow \Sh({\cal C}, J)$ to the functor $f^{\ast} \circ a_{J}\circ y\in {\bf Flat}_{J}({\cal C}, {\cal E})$ (resp. a geometric morphism $f':{\cal E}\rightarrow \Sh(\tilde{{\cal C}}, \tilde{J})$ to the functor $f'^{\ast} \circ a_{\tilde{J}} \circ y\in {\bf Flat}_{\tilde{J}}(\tilde{{\cal C}}, {\cal E})$) (see for example \cite{MM}). From the commutativity of the square above it is then immediate to see that the half of the equivalence given by the composition with the inverse image of $\tau$ corresponds to the obvious restriction functor ${\bf Flat}_{J}({\cal C}, {\cal E})\rightarrow {\bf Flat}_{\tilde{J}}(\tilde{{\cal C}},  {\cal E})$. The other half of the equivalence necessarily induce the functor which sends a $\tilde{J}$-continous flat functor $G:\tilde{\cal C}\rightarrow {\cal E}$ to its extension $\overline{G}:{\cal C}\rightarrow {\cal E}$ to $\cal C$ obtained by putting $\overline{G}(c)=0$ for each $c\in \cal C$ not in $\tilde{\cal C}$; indeed, there is at most one $J$-continuous flat functor $\overline{G}:{\cal C}\rightarrow {\cal E}$ whose restriction to $\tilde{\cal C}$ is a given functor $G:\tilde{\cal C}\rightarrow {\cal E}$ (for $\overline{G}$ to be $J$-continuous, $\overline{G}(c)$ must be equal to $0$ for each object $c\in \cal C$ which is $J$-covered by the empty sieve). 
\end{proofs}
Let us now apply the lemmas above to deduce a model-theoretic characterization of De Morgan (resp. Boolean) toposes among those which arise as localizations of a given presheaf topos $[{\cal C}^{\textrm{op}},\Set]$.\\
Recall that in \cite{OC} we introduced the notion of $J$-homogeneous model of a theory of presheaf type $\mathbb T$ with respect to a Grothendieck topology $J$ on the category $(\textrm{f.p.}{\mathbb T}\textrm{-mod}(\Set))^{\textrm{op}}$. Having this notion in mind, we now introduce the following more specific definition (the notation below being taken from \cite{OC}).\\
Given a Grothendieck topos $\cal E$ with a class of generators $\cal G$, a geometric theory $\mathbb T$ classified by the topos $[{\cal C}^{\textrm{op}}, \Set]$ (with ${\cal C}=(\textrm{f.p.}{\mathbb T}\textrm{-mod}(\Set))^{\textrm{op}}$), and a collection $S$ of arrows in ${\cal C}^{\textrm{op}}$ with common domain, a model $M\in {\mathbb T}\textrm{-mod}(\cal E)$ is said to be $S\textrm{-homogeneous}$ if for each object $E\in \cal G$ and arrow $y:E^{\ast}(\gamma^{\ast}_{\cal E}(i(c)))\rightarrow E^{\ast}(M)$ in ${\mathbb T}\textrm{-mod}({\cal E}\slash E)$ there exists an epimorphic family $(p_{f}:E_{f}\rightarrow E, f\in S)$ and for each arrow $f:c\rightarrow d$ in $S$ an arrow $u_{f}:E_{f}^{\ast}(\gamma^{\ast}_{\cal E}(i(d)))\rightarrow E_{f}^{\ast}(M)$ in ${\mathbb T}\textrm{-mod}({\cal E}\slash E)$ such that $p_{f}^{\ast}(y)=u_{f}\circ E_{f}^{\ast}(\gamma^{\ast}_{\cal E}(i(f)))$.\\
The following results hold.   

\begin{theorem}\label{teom}
Let $\mathbb T$ be a geometric theory classified by the topos $[{\cal C}^{\textrm{op}}, \Set]$ (with ${\cal C}=(\textrm{f.p.}{\mathbb T}\textrm{-mod}(\Set))^{\textrm{op}}$) and $\mathbb T'$ a geometric theory classified by a topos ${\cal E}=\Sh({\cal C}, J)$ together with a full and faithful indexed functor $i:\underline{\mathbb T'\textrm{-mod}}\hookrightarrow \underline{\mathbb T\textrm{-mod}}$ which is induced via the universal property of the classifying toposes by the inclusion $\Sh({\cal C}, J)\hookrightarrow [{\cal C}^{\textrm{op}}, \Set]$. If $m$ is the De Morgan topology on the topos $\cal E$ then $\sh_{m}(\cal E)$ classifies the $\mathbb T'$-models which are $S$-homogeneous (as $\mathbb T$-models via $i$) for each $M_{\tilde{\cal C}}$-covering sieve $S$. In particular, $\mathbb T'$ is a De Morgan theory if and only if every $\mathbb T'$-model (in any Grothendieck topos) is (as a $\mathbb T$-model via $i$) $S$-homogeneous for each $M_{\tilde{\cal C}}$-covering sieve $S$.
\end{theorem}
\begin{proofs}
We have $\sh_{m}({\cal E})=\sh_{m}([\tilde{{\cal C}}^{\textrm{op}}, \Set])\cap {\cal E}$ by Proposition \ref{propDeMorgan}(ii). From this it is clear that $\sh_{m}(\cal E)$ classifies the flat functors on $\tilde{\cal C}$ which are $\tilde{J}$-continuous and $M_{\tilde{\cal C}}$-continuous (where $M_{\tilde{\cal C}}$ is the De Morgan topology on the category $\tilde{\cal C}$), equivalently (by Lemma \ref{lemma2}) the $J$-continuous flat functors on $\cal C$ which send $M_{\tilde{\cal C}}$-covering sieves to epimorphic families. The thesis then follows from Theorems 4.6-4.8 \cite{OC}. 
The last part of the theorem follows from the first part together with Proposition \ref{propDeMorgan}(i) and Lemma \ref{lemma1}.       
\end{proofs}

\begin{theorem}\label{teob}
Let $\mathbb T$ be a geometric theory classified by the topos $[{\cal C}^{\textrm{op}}, \Set]$ (with ${\cal C}=(\textrm{f.p.}{\mathbb T}\textrm{-mod}(\Set))^{\textrm{op}}$) and $\mathbb T'$ a geometric theory classified by a topos ${\cal E}=\Sh({\cal C}, J)$ together with a full and faithful indexed functor $i:\underline{\mathbb T'\textrm{-mod}}\hookrightarrow \underline{\mathbb T\textrm{-mod}}$ which is induced via the universal property of the classifying toposes by the inclusion $\Sh({\cal C}, J)\hookrightarrow [{\cal C}^{\textrm{op}}, \Set]$. Then $\sh_{\neg\neg}(\cal E)$ classifies the $\mathbb T'$-models which are (as $\mathbb T$-models via $i$) $S$-homogeneous for each stably non-empty sieve $S$ in $\tilde{\cal C}$. In particular, $\mathbb T'$ is a Boolean theory if and only if the $\mathbb T'$-models (in any Grothendieck topos) are (identified by $i$ with) the $\mathbb T$-models which are $S$-homogeneous for every stably non-empty sieve $S$ in $\tilde{\cal C}$.
\end{theorem}
\begin{proofs}
This is similar to the proof of Theorem \ref{teom}; we omit the details.
\end{proofs}
It is natural at this point to introduce the following notions.\\
Given a theory $\mathbb T$ classified by a presheaf topos $[{\cal C}^{\textrm{op}}, \Set]$ (with ${\cal C}=(\textrm{f.p.}{\mathbb T}\textrm{-mod}(\Set))^{\textrm{op}}$), we define the theory of De Morgan $\mathbb T$-models to be (the Morita-equivalence class of) the theory of flat functors on $\cal C$ which send $M_{\tilde{\cal C}}$-covering sieves to epimorphic families. Similarly, we define the theory of Boolean $\mathbb T$-models to be (the Morita-equivalence class of) the theory of flat functors on $\cal C$ which send stably non-empty sieves in $\tilde{\cal C}$ to epimorphic families; notice that in case $\cal C$ satisfies the right Ore condition the theory of Boolean $\mathbb T$-models is the same thing as the theory of homogeneous $\mathbb T$-models introduced in \cite{OC2}.\\
So Theorems \ref{teom} and \ref{teob} imply that, up to Morita equivalence, the De Morgan theories are exactly the ``quotients'' of the theories of De Morgan models, while the Boolean theories are precisely the theories of Boolean models.\\
Below we see that the theories of De Morgan and Boolean $\mathbb T$-models admit (up to Morita-equivalence) natural axiomatizations by geometric sequents in the signature of $\mathbb T$.  

\subsection{Axiomatizations}
Let us start with some preliminary observations on the operation of pseudocomplementation in the Heyting algebras of the form $\Sub_{{\cal C}_{\mathbb T}}(\{\vec{x}.\top\})$, where ${\cal C}_{\mathbb T}$ is the syntactic category of a geometric theory $\mathbb T$ over a signature $\Sigma$.\\
We will use the following terminology. A geometric formula $\phi(\vec{x})$ over $\Sigma$ is said to be consistent (with respect to $\mathbb T$) if the sequent $\phi(\vec{x}) \:\vdash_{\vec{x}}\: \bot$ is not provable in $\mathbb T$. Two geometric formulas $\phi(\vec{x})$ and $\psi(\vec{x})$ over $\Sigma$ in the same context are said to be consistent with each other (with respect to $\mathbb T$) if their conjunction is consistent; otherwise, they are said to be inconsistent with each other. $\psi(\vec{x})$ is said to be stably consinstent with $\phi(\vec{x})$ if $\chi(\vec{x})\wedge \phi(\vec{x})$ is consistent for each consistent formula $\chi(\vec{x})$ in the same context which $\mathbb T$-provably implies $\psi(\vec{x})$; $\psi(\vec{x})$ is said to be stably consinstent if it is stably consistent with $\top(\vec{x})$. From now on we will freely identify a geometric formula $\phi(\vec{x})$ over $\Sigma$ with the corresponding monomorphism $\{\vec{x}.\phi(\vec{x})\}\rightarrow \{\vec{x}.\top\}$ in ${\cal C}_{\mathbb T}$.\\
Let us put, for $\phi(\vec{x})\in \Sub_{{\cal C}_{\mathbb T}}(\{\vec{x}.\top\})$, $\textrm{Cons}(\phi(\vec{x}))=\{\psi(\vec{x})\in \Sub_{{\cal C}_{\mathbb T}}(\{\vec{x}.\top\}) \textrm{ | } \psi(\vec{x}) \textrm{ and } \phi(\vec{x}) \textrm{ are consistent}\}$ (and $\textrm{Incons}(\phi(\vec{x})) \textrm{ equal to the complement } \Sub_{{\cal C}_{\mathbb T}}(\{\vec{x}.\top\}) \setminus \textrm{Cons}(\phi(\vec{x}))$). Note that this assignment actually defines a functor $\textrm{Cons}:\Sub_{{\cal C}_{\mathbb T}}(\{\vec{x}.\top\}) \rightarrow {\mathscr{P}}(\Sub_{{\cal C}_{\mathbb T}}(\{\vec{x}.\top\}))$, where ${\mathscr{P}}(\Sub_{{\cal C}_{\mathbb T}}(\{\vec{x}.\top\}))$ is regarded as a poset category with respect to the inclusion. This is in fact an instance of a more general construction, which we describe now.\\
Given a (complete) Heyting algebra $H$, we can define a functor $\textrm{Cons}:H\rightarrow {\mathscr{P}}(H)$ by $\textrm{Cons}(h)=\{h'\in H \textrm{ | } h'\wedge h\neq 0\}$.\\
We can rephrase various concepts involving the operation of pseudocomplementation $\neg$ in $H$ in terms of the functor Cons; for instance, we have the following proposition. 
\begin{proposition}\label{propos}
Let $H$ be a complete Heyting algebra. Then\\
(i) For each $h\in H$, $\neg h=0$ if and only if $\textrm{Cons}(h)=H\setminus \{0\}$\\
(ii) For each $h,h'\in H$, $h'\leq \neg\neg h$ if and only if $\textrm{Cons}(h')\subseteq \textrm{Cons}(h)$, if and only if  $h'$ is stably joint with $h$ (i.e. $a\wedge h\neq 0$ for each $a\neq 0$ such that $a\leq h'$).\\
(iii) For each $h\in H$, $\neg\neg h=h$ if and only if for every $h'\in H$, $\textrm{Cons}(h')\subseteq\textrm{Cons}(h)$ implies $h'\leq h$.\\
(iv) $H$ is a Boolean algebra if and only if the functor Cons is conservative.
\end{proposition}
\begin{proofs}
(i) This is immediate from the fact that for each $h\in H$, $\neg h=\mathbin{\mathop{\textrm{\huge $\vee$}}\limits_{h'\wedge h=0}h'}$.\\
(ii) For each $h,h'\in H$, $h'\leq \neg\neg h$ if and only if $h'\wedge \neg h=0$, if and only if  $h'\wedge \mathbin{\mathop{\textrm{\huge $\vee$}}\limits_{a\wedge h=0}a}=\mathbin{\mathop{\textrm{\huge $\vee$}}\limits_{a\wedge h=0}a\wedge h'}=0$, if and only if $\textrm{Cons}(h')\subseteq \textrm{Cons}(h)$. The last equivalence in (ii) is obvious.\\ 
(iii) For any $h\in H$, $\neg\neg h=h$ if and only if $\neg\neg h\leq h$, if and only if for each $h'\in H$ $h'\leq \neg\neg h$ implies $h'\leq h$, if and only if $\textrm{Cons}(h')\subseteq\textrm{Cons}(h)$ implies $h'\leq h$, where the last equivalence follows from (ii).\\ 
(iv) $H$ is a Boolean algebra if and only for each $h\in H$ $\neg h=0$ implies $h=1$ (one direction is obvious, while for the other one it suffices to observe that $h\vee \neg h=1$ for each $h\in H$ as $\neg(h\vee \neg h)=\neg h\wedge \neg\neg h=0$), if and only if for each $h\in H$ $\neg\neg h=h$; then our thesis follows at once from (i) and (iii).    
\end{proofs}
Given a complete Heyting algebra $H$, we note that $H$ is a De Morgan algebra if and only if for each pair of elements $h, h'\in H$, $h\wedge h'=0$ implies $\neg h\vee \neg h'=1$ that is $\mathbin{\mathop{\textrm{\huge $\vee$}}\limits_{a\wedge h=0 \textrm{ or } a\wedge h'=0 }a}=1$; as we have observed above, $H$ is a Boolean algebra if and only if for each $h\in H$ $\neg h=0$ implies $h=1$.\\
These characterizations are the ingredients of our axiomatizations.
\newpage
\begin{theorem}\label{axiom1}
Let $\mathbb T$ be a geometric theory over a signature $\Sigma$ classified by a Grothendieck topos $\cal E$. Then the theory $\mathbb T'$ obtained by adding to the axioms of $\mathbb T$ all the geometric sequents of the form $\top \:\vdash_{\vec{x}}\: \mathbin{\mathop{\textrm{\huge $\vee$}}\limits_{\psi(\vec{x})\in \textrm{Incons}(\phi(\vec{x}))\cup\textrm{Incons}(\phi'(\vec{x}))}\psi(\vec{x})}$, where $\phi(\vec{x})$ and $\phi'(\vec{x})$ are geometric formulas in the same context which are inconsistent with each other with respect to $\mathbb T$, is classified by the topos $\sh_{m}(\cal E)$ (where $m$ is the De Morgan topology on $\cal E$). 
\end{theorem}
\begin{proofs}
Let us represent $\cal E$ as $\Sh({\cal C}_{\mathbb T}, J_{\mathbb T})$. Recall that we have an equivalence of categories ${\mathbb T}\textrm{-mod}({\cal E})\simeq {\bf Flat}_{J_{\mathbb T}}({\cal C}_{\mathbb T}, {\cal E})$ (natural in ${\cal E}\in \mathfrak{BTop}$) which sends each model $M\in {\mathbb T}\textrm{-mod}({\cal E})$ the functor $F:{\cal C}_{\mathbb T}\rightarrow {\cal E}$ assigning to a formula $\phi(\vec{x})$ its interpretation $[[\phi(\vec{x})]]_{M}$ in $M$. 
As we have observed in the proof of Theorem \ref{teom}, $\sh_{m}(\cal E)$ classifies the $J_{\mathbb T}$-continuous flat functors on ${\cal C}_{\mathbb T}$ which send $M_{\tilde{{\cal C}_{\mathbb T}}}$-covering sieves to epimorphic families; it remains to show that, via the equivalence above, these functors correspond precisely to the $\mathbb T$-models $M$ such that $[[\neg\phi(\vec{x}) \vee \neg\neg \phi(\vec{x}) ]]_{M}=[[\top(\vec{x})]]_{M}$ for each geometric formula $\phi(\vec{x})$ over $\Sigma$ (equivalently, $[[\neg\phi(\vec{x}) \vee \neg \psi(\vec{x}) ]]_{M}=[[\top(\vec{x})]]_{M}$ for each pair $\phi(\vec{x})$ and $\psi(\vec{x})$ of geometric formulas over $\Sigma$ in the same context which are inconsistent with each other with respect to $\mathbb T$). By our results in section $2$ and Lemma \ref{lemma1}, we have that the $J_{\mathbb T}$-continuous flat functors on ${\cal C}_{\mathbb T}$ which send $M_{\tilde{{\cal C}_{\mathbb T}}}$-covering sieves to epimorphic families are exactly those which send every family of arrows of the form $M_{\phi(\vec{x})}=\{\psi(\vec{x}) \textrm{ | ($\psi(\vec{x})\in \textrm{Incons}(\phi(\vec{x}))$ or ($\psi(\vec{x})$ is stably consistent with $\phi(\vec{x})$)}\}$
 (for a geometric formula $\phi(\vec{x})$ over $\Sigma$) to an epimorphic family; by Proposition \ref{propos}(ii), this is clearly equivalent to saying that the corresponding models $M$ satisfy $[[\neg\phi(\vec{x}) \vee \neg\neg \phi(\vec{x}) ]]_{M}=[[\top(\vec{x})]]_{M}$ for each geometric formula $\phi(\vec{x})$ over $\Sigma$.     
\end{proofs}
The theory $\mathbb T'$ defined in the theorem above will be called the \emph{DeMorganization} of the theory $\mathbb T$.

\begin{theorem}\label{axiom2}
Let $\mathbb T$ be a geometric theory over a signature $\Sigma$ classified by a Grothendieck topos $\cal E$. Then the theory $\mathbb T'$ obtained by adding to the axioms of $\mathbb T$ all the geometric sequents of the form $\top \:\vdash_{\vec{x}}\: \phi(\vec{x})$, where $\phi(\vec{x})$ is a geometric formula over $\Sigma$ which is stably consistent with respect to $\mathbb T$, is classified by the topos $\sh_{\neg\neg}(\cal E)$. 
\end{theorem}
\begin{proofs}
The proof proceeds analogously to that of Theorem \ref{axiom1}, by using Theorem \ref{teob} and the fact that the $J_{\mathbb T}$-continuous flat functors on ${\cal C}_{\mathbb T}$ which send stably non-empty sieves in $\tilde{{\cal C}_{\mathbb T}}$ to epimorphic families are exactly those which send each sieve in $\tilde{{\cal C}_{\mathbb T}}$ generated by a stably consistent formula (cfr. Section 2) to an epimorphic family.
\end{proofs}
The theory $\mathbb T'$ defined in the theorem above will be called the \emph{Booleanization} of the theory $\mathbb T$.

From the theorems above we can deduce the following corollaries.
\begin{corollary}
Let $\mathbb T$ be a theory classified by a topos $[{\cal C}^{\textrm{op}}, \Set]$ (with ${\cal C}=(\textrm{f.p.}{\mathbb T}\textrm{-mod}(\Set))^{\textrm{op}}$). Then the DeMorganization of $\mathbb T$ axiomatizes the $\mathbb T$-models which are $S$-homogeneous for every $M_{\cal C}$-covering sieve $S$ (where $M_{\cal C}$ is the De Morgan topology on the category ${\cal C}$). 
\end{corollary}
\begin{proofs}
This immediately follows from Theorems \ref{teom} and \ref{axiom1} by using Theorems 4.6-4.8 \cite{OC}.
\end{proofs}

\begin{corollary}
Let $\mathbb T$ be a theory classified by a topos $[{\cal C}^{\textrm{op}}, \Set]$ (with ${\cal C}=(\textrm{f.p.}{\mathbb T}\textrm{-mod}(\Set))^{\textrm{op}}$). Then the Booleanization of $\mathbb T$ axiomatizes the $\mathbb T$-models which are $S$-homogeneous for every stably non-empty sieve $S$ in the category ${\cal C}$. In particular, if $\cal C$ satisfies the right Ore condition, the Booleanization of $\mathbb T$ classifies the homogeneous $\mathbb T$-models.  
\end{corollary}
\begin{proofs}
This immediately follows from Theorems \ref{teob} and \ref{axiom2} by using Theorems 4.6-4.8 \cite{OC}.
\end{proofs}

\section{Examples}

In this section we provide some examples of theories which are De Morgan and theories which are not, focusing our attention on coherent theories (recall that we have proved that every regular theory is De Morgan, cfr. Theorem \ref{teosyintreg} above).\\
The first theory we consider is the theory of dense linear orders. As remarked in \cite{blasce}, this theory is not Boolean; however it is a De Morgan theory, as stated in the following proposition.

\begin{proposition}
The theory of dense linear orders is a De Morgan theory.
\end{proposition}
 
\begin{proofs}
By an obvious variation of the arguments in Example D3.4.11 \cite{El2}, the classifying topos for the theory of dense linear orders can be represented as the topos $\Sh({{\bf Ord}^{\textrm{op}}_{fm}}, J)$ of sheaves on the opposite of the category ${\bf Ord}_{fm}$ of finite ordinals and order-preserving injections between them with respect to a Grothendieck topology $J$ on ${{\bf Ord}^{\textrm{op}}_{fm}}$ with no empty covering sieves. Note that the category ${\bf Ord}_{fm}^{\textrm{op}}$ satisfies the right Ore condition; our thesis then follows from Corollary \ref{corfond}.
\end{proofs} 
Next, let us consider a couple of theories which arise as ``quotients'' of the algebraic theory of rings. 
 
\begin{proposition}
The theory of local rings is not a De Morgan theory.
\end{proposition}
 
\begin{proofs}
It is well-known that the coherent theory of local rings is classified by the Zariski topos $\cal Z$, that is the topos $\Sh({{\bf Rng}_{f.g.}^{\textrm{op}}}, J)$ of sheaves on the opposite of the category ${\cal C}={\bf Rng}_{f.g.}$ of finitely generated rings with respect to the topology $J$ on ${\bf Rng}_{f.g.}^{\textrm{op}}$ defined as follows: given a cosieve $S$ in ${\bf Rng}_{f.g.}$ on an object $A$, $S\in J(A)$ if and only if $S$ contains a finite family $\{\xi_{i}:A\rightarrow A[{s_{i}}^{-1}] \textrm{ | } 1\leq i\leq n\}$ of canonical inclusions $\xi_{i}:A\rightarrow A[{s_{i}}^{-1}]$ in ${\bf Rng}_{f.g.}$ where $\{s_{1},\ldots, s_{n}\}$ is any set of elements of $A$ which is not contained in any proper ideal of $A$. Obviously, the only object of ${\bf Rng}_{f.g.}$ which is $J$-covered by the empty cosieve is the zero ring. In order to apply our criterion `$\Sh({\cal C}, J)$ is De Morgan if and only if $M_{\tilde{{\cal C}}}\leq \tilde{J}$' (established in section $2$ above) to decide whether $\cal Z$ is De Morgan or not, let us consider the reduced site $(\tilde{{\cal C}}, \tilde{J})$. Notice that for $A\in \tilde{{\cal C}}$, the canonical inclusion $\xi_{s}:A\rightarrow A[s^{-1}]$ lies in $\tilde{{\cal C}}$ (i.e. it is not the zero map) if and only if $s$ is not nilpotent. Then we have: $S\in \tilde{J}(A)$ if and only if $S$ contains a finite family $\{\xi_{s_{i}}:A\rightarrow A[{s_{i}}^{-1}] \textrm{ | } 1\leq i\leq n\}$ of canonical inclusions $\xi_{s_{i}}:A\rightarrow A[{s_{i}}^{-1}]$ in ${\bf Rng}_{f.g.}$ where $\{s_{1},\ldots, s_{n}\}$ is any set of non-nilpotent elements of $A$ which is not contained in any proper ideal of $A$. Indeed, by definition of induced topology, $S\in \tilde{J}(A)$ if and only if there exists elements $s_{1}, \ldots, s_{n}$ and $t_{1}, \ldots, t_{m}$ such that the set $\{s_{1}, \ldots, s_{n}, t_{1}, \ldots, t_{m} \}$ is not contained in any proper ideal of $A$, all the $s_{i}$ are non-nilpotent, all the $t_{j}$ are nilpotent and $S$ contains the family $\{\xi_{s_{1}},\ldots \xi_{s_{1}}, \xi_{t_{1}},\ldots \xi_{t_{m}} \}$; but $(s_{1}, \ldots, s_{n}, t_{1}, \ldots, t_{m})=1$ implies $(s_{1}^{r}, \ldots, s_{n}^{r}, t_{1}^{r}, \ldots, t_{m}^{r})=1$ for any natural number $r$, in particular for an $r$ such that $t_{j}^{r}=0$ for each $1\leq j\leq m$ and hence $(s_{1}, \ldots, s_{n})=1$.\\
For an object $A\in \tilde{{\cal C}}$ and a non-nilpotent element $a\in A$, let us denote by $S^{A}_{a}$ the cosieve in $\tilde{{\cal C}}$ generated by the arrow $\xi_{a}:A\rightarrow A[a^{-1}]$; note that $S^{A}_{a}$ can be identified with the collection of arrows $f$ in $\tilde{{\cal C}}$ with domain $A$ such that $f(a)$ is invertible. Let us now observe some facts about these sieves.

\begin{lemma}\label{lemmazariski}
With the above notation, for any object $A\in \tilde{{\cal C}}$ and non-nilpotent element $a\in A$ the following facts hold:\\
(i) For any arrow $f:A\rightarrow B$ in $\tilde{{\cal C}}$ with domain $A$,\\
\[
\begin{array}{rcll}
f^{\ast}(S^{A}_{a}) & = & \emptyset & \textrm{if $f(a)$ is nilpotent}\\
                & = & S^{B}_{f(a)} & \textrm{if $f(a)$ is not nilpotent;}\\

\end{array}
\]
(ii) $S^{A}_{a}$ is a $\tilde{J}$-closed cosieve;\\
(iii) $S^{A}_{a}$ is stably non-empty if and only if it is the maximal cosieve on $A$ i.e. $a$ is invertible.\\
\end{lemma}

\begin{proofs}
(i) This is immediate from the equalities $f^{\ast}(S^{A}_{a})=\{g:B\rightarrow cod(g) \textrm{ | } g\circ f\in S^{A}_{a}\}=\{g:B\rightarrow cod(g) \textrm{ | $g(f(a))$ is invertible}\}$.\\
(ii) Suppose that for an arrow $f:A\rightarrow B$ in $\tilde{{\cal C}}$ with domain $A$ we have $f^{\ast}(S^{A}_{a})\in \tilde{J}(B)$; then (by the description of $\tilde{J}$ given above) there exist a finite number of non-nilpotent elements $s_{1}, \ldots, s_{n}\in B$ such that $(s_{1}, \ldots, s_{n})=1$ and $\xi_{s_{i}}$ belongs to $f^{\ast}(S^{A}_{a})$ for each $1\leq i\leq n$. So we have (by the calculation above) that $\xi_{s_{i}}(f(a))$ is invertible for each $1\leq i\leq n$; this in turn means that there exists for each $1\leq i\leq n$ an element $c_{i}\in B$ and a natural number $n_{i}$ such that $f(a)c_{i}=s_{i}^{n_{i}}$. Now, by taking $p=max\{n_{i}\}$ we have $(s_{1}^{p}, \ldots, s_{n}^{p})=1$ and hence the existence of an element $c$ such that $f(a)c=1$, that is $f(a)$ invertible in $B$ i.e. $f$ belongs to $S^{A}_{a}$.\\
(iii) One direction is obvious; let us prove the other one. If $a$ is not invertible then $(a)$ is a proper ideal of $A$ and hence the quotient $A/(a)$ belongs to $\tilde{{\cal C}}$. But if $\pi:A\rightarrow A/(a)$ is the natural projection map we have (by part (i) of the lemma) $\pi^{\ast}(S^{A}_{a})=\emptyset$, hence $S^{A}_{a}$ is not stably non-empty.           
\end{proofs}
Now by Lemma \ref{lemmazariski} we have that, for a given object $A\in \tilde{{\cal C}}$ and non-nilpotent element $a\in A$, $M_{(\xi_{a})}=\{f:A\rightarrow cod(f) \textrm{ | $f(a)$ is nilpotent or invertible}\}$ (with the notation of section 2 above). Suppose now that $A$ is an integral domain and $M_{(\xi_{a})}$ is $\tilde{J}$-covering. Then there exists a finite number of non-nilpotent elements $s_{1}, \ldots, s_{n}$ of $A$ such that $(s_{1}, \ldots, s_{n})=1$ and $\xi_{s_{i}}$ belongs to $M_{(\xi_{a})}$ for each $1\leq i\leq n$. If for some $i$ $\xi_{s_{i}}(a)$ were nilpotent then, $A$ being an integral domain, $s_{i}$ would be nilpotent, contradicting our assumption. So we deduce that $\xi_{s_{i}}(a)$ is invertible for each $1\leq i\leq n$, from which it follows (as in the proof of part (i) of Lemma \ref{lemmazariski}) that $a$ is invertible. So we have proved that if $A\in \tilde{{\cal C}}$ is an integral domain and $a\in A$ is a non-nilpotent element of $A$ then $M_{(\xi_{a})}$ is $\tilde{J}$-covering if and only $a$ is invertible. Any instance of an integral domain $A$ together with a non-invertible element $a\in A$ (for example $A$ equal to the ring of integers $\mathbb Z$ and $a=2$) then proves that $\cal Z$ is not De Morgan.   
\end{proofs}

\begin{proposition}
The theory of integral domains is not a De Morgan theory.
\end{proposition}
 
\begin{proofs}
The theory $\mathbb T$ of integral domains is obtained from the algebraic theory of rings by adding the following axioms:
\[
0=1 \: {\vdash}_{\{\}}\: \bot;
\]
\[
x \cdot y =0 \ \: \vdash_{x,y}\: (x=0)\vee (y=0).
\]
By using Proposition D3.1.10 \cite{El2} we immediately obtain the following representation for the classifying topos ${\Set}[{\mathbb T}]$ of $\mathbb T$: ${\Set}[{\mathbb T}]\simeq \Sh({{\bf Rng}_{f.g.}^{\textrm{op}}}, J)$, where $J$ is the smallest Grothendieck topology on ${\cal C}={\bf Rng}_{f.g.}^{\textrm{op}}$ such that the empty sieve on the zero ring and the cosieve in ${\bf Rng}_{f.g.}$ on ${\mathbb Z}[x,y]/(x\cdot y)$ generated by the canonical projections ${\mathbb Z}[x,y]/(x\cdot y)\rightarrow {\mathbb Z}[x,y]/(x)$ and ${\mathbb Z}[x,y]/(x\cdot y)\rightarrow {\mathbb Z}[x,y]/(y)$ are $J$-covering.
The following result (which was motivated by the arguments p. 111-112 \cite{MM}) will be useful for giving an explicit description of the topology $J$.
\begin{lemma}\label{grot}
Let $\cal C$ be a category and $K$ be a function which assigns to each object $c\in {\cal C}$ a collection $K(c)$ of sieves in $\cal C$ on $c$. Then $K$ is a Grothendieck topology if and only if it satisfies the following properties:\\
(i) the maximal sieve $M(c)$ belongs to $K(c)$;\\
(ii) for each pair of sieves $S$ and $T$ on $c$ such that $T\in K(c)$ and $S\supseteq T$, $S\in K(c)$;\\
(iii) if $R\in K(c)$ then for any arrow $g:d\rightarrow c$ there exists a sieve $S\in K(c)$ such that for each arrow $f$ in $S$, $g\circ f\in R$;\\
(iv) if $\{f_{i}:c_{i}\rightarrow c \textrm{ | } i\in I \}\in K(c)$ and for each $i\in I$ we have a sieve $\{g_{ij}:d_{ij}\rightarrow c_{i} \textrm{ | } j\in I_{i} \}\in K(c_{i})$, then the family of composites $\{f_{i}\circ g_{ij}:d_{ij}\rightarrow c \textrm{ | } i\in I, j\in I_{i} \}$ belongs to $K(c)$. 
\end{lemma}
\begin{proofs}
In one direction, let us suppose that $K$ is a Grothendieck topology. Properties (i) and (ii) are well-known to hold. Property (iii) holds as we can clearly take as cover $S$ satisfying the condition the pullback of the sieve $R$ along the arrow $g:d\rightarrow c$, which is $K$-covering by the stability axiom for Grothendieck topologies. Property (iv) easily follows from the transitivity axiom for Grothendieck topologies; indeed, the sieve $R:=\{f_{i}\circ g_{ij}:d_{ij}\rightarrow c \textrm{ | } i\in I, j\in I_{i} \}$ satisfies the following property with respect to the sieve $S:=\{f_{i}:c_{i}\rightarrow c \textrm{ | } i\in I \}\in K(c)$: for all arrows $h$ in $S$, $h^{\ast}(R)$ is $K$-covering.\\
Conversely, let us suppose that $K$ satisfies properties (i), (ii), (iii) and (iv). To prove that $K$ satisfies the stability axiom for Grothendieck topologies we observe that if $R\in K(c)$ and $g:d\rightarrow c$ is an arrow with codomain $c$, then $h^{\ast}(R)$ contains the sieve $S$ given by property (iii) and hence is $K$-covering by property (ii). It remains to verify that $K$ satisfies the transitivity axiom for Grothendieck topologies. Given a sieve $R$ on $c$ and a sieve $S\in K(c)$ such that for all arrows $h$ in $S$, $h^{\ast}(R)$ is $K$-covering, we want to prove that $R$ is $K$-covering. This follows from property (ii) as $R$ contains ``the composite`` of the sieve $S$ with the sieves of the form $h^{\ast}(R)$ for $h$ in $S$.    
\end{proofs}
Notice that, in case $\cal C$ has pullbacks, property (iii) in the lemma may be replaced by the following condition: if $\{f_{i}:c_{i}\rightarrow c \textrm{ | } i\in I \}\in K(c)$ then for any arrow $g:d\rightarrow c$ the sieve generated by the family of pullbacks $\{ \textrm{p.b.($f_{i}$, $g$)} \rightarrow d \textrm{ | } i\in I \}$ belongs to $K(d)$.\\

Now, as a consequence of Lemma \ref{grot}, it is immediate to see that $J$ is the topology defined as follows: given a cosieve $S$ in ${\bf Rng}_{f.g.}$ on an object $A\in {\bf Rng}_{f.g.}$, $S\in J(A)$ if and only if either $A$ is the zero ring and $S$ is the empty sieve on it or $S$ contains a finite family $\{\pi_{a_{i}}:A\rightarrow A/(a_{i}) \textrm{ | } 1\leq i\leq n\}$ of canonical projections $\pi_{a_{i}}:A\rightarrow A/(a_{i})$ in ${\bf Rng}_{f.g.}$ where $\{a_{1},\ldots,a_{n}\}$ is any set of elements of $A$ such that $a_{1}\cdot \ldots \cdot a_{n}=0$.\\
Let us now observe the following fact.
\begin{lemma}\label{countroc}
Let $\cal C$ be a category and $J$ be a Grothendieck topology on it with no empty covering sieves. If there exists an object $c\in {\cal C}$ such that $J(c)=\{M(c)\}$ and two arrows $f$ and $g$ with codomain $c$ such that $f^{\ast}((g))=\emptyset$ then $\Sh({\cal C}, J)$ is not a De Morgan topos. 
\end{lemma}
\begin{proofs}
By using the notation of section $2$ above, we have that $M_{(g)}\in J(c)$ implies $(g)$ empty or stably non-empty; as none of the two alternatives hold in our case (in view of our hypotheses), we conclude that $M_{(g)}\notin J(c)$ and hence, by Theorem \ref{criterion}, $\Sh({\cal C}, J)$ is not a De Morgan topos.  
\end{proofs}
Lemma \ref{countroc} provides us with a counterexample to $\Set[{\mathbb T}]$ being De Morgan. Indeed, we note that if $A\in {\bf Rng}_{f.g.}$ is an integral domain then $\tilde{J}(A)=\{M(A)\}$; so by taking $A={\mathbb Z}$, and $f$ and $g$ to be respectively the canonical projections ${\mathbb Z}\rightarrow {\mathbb Z}/2{\mathbb Z}$ and ${\mathbb Z}\rightarrow {\mathbb Z}/3{\mathbb Z}$ in ${\bf Rng}_{f.g.}$ we have that the reduced site $(\tilde{\cal C}, \tilde{J})$ satisfies the hypotheses of the lemma.
\end{proofs} 

\vspace{3 mm}
An analysis of the theory of fields in relation to De Morgan's law will follow shortly in a separate paper. 
\vspace{7 mm}

{\bf Acknowledgements:} I am very grateful to my Ph.D. supervisor Peter Johnstone for suggesting me to investigate De Morgan classifying toposes and for his support during the preparation of this work. Thanks also to Martin Hyland for useful discussions.

\end{document}